\documentclass[3pt]{article}

\usepackage{graphicx}\usepackage{amssymb}
\usepackage{amsmath}\usepackage{amsfonts}
\numberwithin{equation}{section}

\begin{document}

\title{The lower bound property of the Morley element
eigenvalues }
\author{ { Yidu Yang, Hao Li, Hai Bi} \\
{\small School of Mathematics and Computer Science, }\\{\small
Guizhou Normal University,  Guiyang,  $550001$,  China}\\{\small
ydyang@gznu.edu.cn, lihao365@126.com, bihaimath@gznu.edu.cn}
}\date{~} \pagestyle{plain} \textwidth 145mm \textheight 215mm
\topmargin 0pt \maketitle

\indent{\bf\small Abstract~:~} {\small In this paper,
 we prove that the Morley element
eigenvalues approximate the exact ones from below on regular meshes,
including adaptive local refined meshes, for the fourth-order
elliptic eigenvalue problems with the clamped boundary condition in
any dimension. And we implement the adaptive computation to obtain
lower bounds of the Morley element eigenvalues for the vibration
problem of clamped plate under tension.
}\\
\indent{\bf\small Keywords~:~\scriptsize} {\small fourth-order
elliptic eigenvalue problems, Morley elements, regular mesh, lower
bounds of eigenvalues, adaptive computation.
\\\indent{\bf\small 1991 MSC~:~\scriptsize}    code 65N25,65N30}

\section{Introduction}

\indent The Morley element is a non-conforming triangle element
proposed by Morley \cite{morley} in 1968 for plate bending problems.
Also, the Morley element was extended to arbitrarily dimensions by
Wang and Xu \cite{wm,wm2}. So, this element is also called the Morley-Wang-Xu element (see \cite{hu1}).\\
\indent  Using the Morley element to obtain the lower bounds for
eigenvalues of fourth-order elliptic eigenvalue problems is a
problem of concern from mathematical and mechanical community. In
1979 Rannacher \cite{rr2} found through numerical computation that
the Morley element can obtain the lower eigenvalue bounds for the
vibration of clamped plate. This discovery is very important in
engineering and mechanics computing. Lin et al. \cite{yang2} first
proved this discovery theoretically. Hu et al. \cite{hu1} extended
the work in \cite{yang2} to $2m$-th order elliptic eigenvalue
problems in arbitrary dimensions, and Lin et al. \cite{lin2} also
developed the work
in \cite{yang2} further.\\
\indent The references \cite{hu1,lin2,yang2} studied the lower bound property in the asymptotic sense.
How to check that whether the mesh size
is small enough or not in practice? In fact, in computation the
Morley eigenvalues will become more and more precise when the mesh
is refined gradually. So it can be concluded that the condition that
the mesh size is small enough is satisfied when the Morley element
eigenvalues  reveal a stable  monotonically increasing tendency.
Thus, it can be deduced that the Morley element gives the lower
eigenvalue bounds. The numerical examples in \cite{carstensen,rr2}
and Section 5 in this paper all support this conclusion.
Hence, it is a very meaningful work to study the asymptotic lower bound. \\
\indent It is noteworthy that Carstensen and Gallistl
\cite{carstensen} studied the guaranteed lower eigenvalue bounds for
the biharmonic equation by a simple post-processing method for the
Morley element eigenvalues. Although the eigenvalues corrected in
\cite{carstensen} may not accurate than the Morley eigenvalues (see the numerical experiments report therein), the work \cite{carstensen} is also very
important and meaningful.\\
\indent A posteriori error estimates and adaptive methods of finite
element approximation are topics attracting more attention from
mathematical and physical fields (see, e.g.,
\cite{ainsworth,ainsworth2,babuska3,becker2,carstensen2,dai1,duran,gedicke,morin,verfurth}
and the references cited therein), and have also been applied to the
Morley element method for plate problems (see, for example,
\cite{bl,gallistl,hu2,yang5}).\\
\indent Based on the above work, this paper further studies the
asymptotic lower bound property of the Morley element eigenvalues on
regular meshes, including adaptive local refined meshes.
The features of this paper are as follows:\\
\indent (1)~For fourth-order elliptic eigenvalue problems with the
clamped boundary condition in any dimension, including the
vibrations of a clamped plate under tension, we prove in the
asymptotic sense that the Morley element eigenvalues approximate the
exact ones from below.\\
\indent (2)~Under the saturation condition $\|u-u_{h}\|_{h}\gtrsim
h^{t_{0}}$, \cite{hu1,yang2} studied the approximation from below
where $t_{0}$ is the singularity exponent of the eigenfunction $u$.
However, this condition is not valid on adaptive meshes with local
refinement. \cite{lin2} discussed the approximation from below on
quasi-uniform meshes and gave the stability condition
$\|u-u_{h}\|_{h}\gtrsim h^{2}$. Developing the work in
\cite{hu1,lin2,yang2} we prove on regular meshes, including adaptive
local refined meshes, that the lower bound property of the Morley
element
eigenvalue.\\
\indent (3)~Thanks to \cite{dai1,yang5}, we get the relationship
between the Morley element eigenvalue approximation and the
associated Morley element boundary value approximation with
$\lambda_{h}u_{h}$ as the right hand term, thus we obtain for the
vibration problem of a clamped plate the reliable and effective a
posteriori error estimators for the Morley element eigenpair which
come from those given by Hu and Shi \cite{hu2} for the plate bending
problem. Shen \cite{shen} also discussed a posteriori error estimators for
the Morley element eigenpair, while we focus in this paper on the reliability
and effectiveness of the a posteriori error estimators on adaptive meshes.
And thus, based on these a posteriori error estimator we implement the adaptive
computation using the package of iFEM \cite{chenl}. The numerical
results validate that the a posteriori error estimators are sharp,
and the lower bound property of the Morley element eigenvalues on
adaptive meshes.\\

\indent In this paper, regarding the basic theory of  finite
elements, we refer to \cite{babuska,brenner,ciarlet,oden,shi2}.\\

\indent Throughout this paper, $C$ denotes a positive constant
independent of mesh size, which may not be the same in different
places. For simplicity, we use the notation $a\lesssim b$ to mean
that $a\leq C b$, and use $a\thickapprox b$ to mean
that $a\lesssim b$ and $b\lesssim a$. Finally, $a\backsimeq b$
abbreviates $a=b$ in the asymptotic sense.

\setcounter{section}{1}\setcounter{equation}{0}

\section{preliminaries.}

\indent Consider the fourth-order elliptic eigenvalue problem:
\begin{eqnarray}\label{2.1}
&&Lu\equiv \sum\limits_{i,j=1}^{n}\frac{\partial^{2}}{\partial
x_{j}\partial x_{i}}(a_{ij} \frac{\partial^{2}u}{\partial
x_{i}\partial x_{j}})
-\sum\limits_{i,j=1}^{n}\frac{\partial}{\partial x_{j}}(b_{ij}
\frac{\partial u}{\partial x_{i}}) +\beta u=\lambda \rho u,~in
\Omega,~~~~~~\\\label{2.2}
 &&~~~~~~~~\frac{\partial u}{\partial
\gamma}=0,~~~u=0,~~~on~ \partial\Omega,
\end{eqnarray}
where $\Omega\subset R^{n}~(n=2,3,\cdots)$ is a polyhedral domain
with boundary $\partial\Omega$, $\frac{\partial u}{\partial\gamma}$
is the outward normal derivative on $\partial\Omega$,
$(a_{ij})_{n\times n}$ and $(b_{ij})_{n\times n}$  are symmetric
matrices, $a_{ij}$ and $b_{ij}$ are appropriate smooth functions,
$\beta, \rho\in L_{\infty}(\Omega)$, $\sum b_{ij}\xi_{i}\xi_{j}\geq
0$, $\forall x\in\Omega$ and $\forall \xi\in R^{n}$, $0\leq \beta$,
$a_{ij}$$(i,j=1,2,\cdots,n)$ and $\rho$ are bounded below by a
positive constant
on $\Omega$.\\
\indent Let $W_{m,p}(\Omega)$ be a Sobolev space with norm
$\|\cdot\|_{m,p}$ and semi-norm $|\cdot|_{m,p}$,
$H^{m}(\Omega)=W_{m,2}(\Omega)$, $H_{0}^{2}(\Omega)=\{v\in
H^{2}(\Omega): v|_{\partial\Omega}=\frac{\partial v}{\partial \gamma}|_{\partial\Omega}=0\}$.\\
\indent The weak form of (\ref{2.1})-(\ref{2.2}) is to seek
$(\lambda,u)\in R\times H_{0}^{2}(\Omega)$ with $\|u\|_{b}=1$ such
that
\begin{eqnarray}\label{2.3}
a(u,v)=\lambda b(u,v),~~~\forall v\in H_{0}^{2}(\Omega),
\end{eqnarray}
where
\begin{eqnarray*}
&&a(u,v)= \int\limits_{\Omega}(\sum\limits_{i,j=1}^{n} a_{ij}
\frac{\partial^{2}u}{\partial x_{i}\partial
x_{j}}\frac{\partial^{2}v}{\partial x_{i}\partial x_{j}}
+\sum\limits_{i,j=1}^{n}b_{ij} \frac{\partial u}{\partial
x_{i}}\frac{\partial v}{\partial x_{j}}
+\beta u v) dx,\\
&&b(u,v)=\int\limits_{\Omega}\rho uv dx,~~~\|u\|_{b}=\sqrt{b(u,u)}.
\end{eqnarray*}
It is obvious that $a(u,v)$ is a symmetric, continuous, and
$H_{0}^{2}(\Omega)$-elliptic bilinear form, and $b(u,v)$ is a
symmetric, continuous and positive definite bilinear form. Let
$\|u\|_{a}=\sqrt{a(u,u)}$. Then the norms $\|u\|_{a}$, $\|u\|_{2,2}$, and $|u|_{2,2}$ are equivalent, and $\|u\|_{b}$ is equivalent to $\|u\|_{0,2}$.\\
\indent We assume $\pi_{h}=\{\kappa\}$ is a regular simplex
partition of $\Omega$ and satisfies
$\overline{\Omega}=\bigcup\overline{\kappa}$ (see \cite{ciarlet}).
Let $h_{\kappa}$ be the diameter of $\kappa$, and
$h=\max\{h_{\kappa}: \kappa\in \pi_{h}\}$ be the mesh size of
$\pi_{h}$. Let $\varepsilon_{h}=\{F\}$ denotes the set of faces
(($n-1$)-simplexes ) of $\pi_{h}$, and let $\varepsilon_{h}'=\{l\}$
denotes the set of faces ($n-2$)-simplexes of $\pi_{h}$. When $n=2$,
$l=z$ is a vertex of $\kappa$, and
\begin{eqnarray*}
\frac{1}{meas(l)}\int\limits_{l}v=v(z).
\end{eqnarray*}
Let $\pi_{h}(\kappa)$ denotes the set of all elements sharing common face with the element $\kappa$.\\
\indent Let $\kappa_{+}$ and $\kappa_{-}$ be any two n-simplex with
a face $F$ in common such that the unit outward normal to
$\kappa_{-}$ at $F$ corresponds to $\gamma_{F}$. We denote the jump
of $v$ across the face $F$ by
$$[v]=(v|_{\kappa_{+}}-v|_{\kappa_{-}})|_{F}.$$
And the jump on boundary faces is simply given by the trace of the function on each face.\\
\indent In \cite{wm}, the Morley-Wang-Xu element space is defined by
\begin{eqnarray*}
&&S^{h}=\{v\in L_{2}(\Omega):v\mid_{\kappa} \in
P_{2}(\kappa),\forall
\kappa\in \pi_{h},\\
&&~~~~~~\int\limits_{F}[\nabla v\cdot \gamma_{F}]=0~~\forall F\in
\varepsilon_{h}, ~~~\int\limits_{l}[v]=0~~\forall l\in
\varepsilon_{h}' \},
\end{eqnarray*}
where $P_{2}(\kappa)$ denotes the space of polynomials of degree less than or equal to 2 on $\kappa$.\\
\indent The Morley-Wang-Xu element space $S^{h}\subset L_{2}(\Omega)$ and $S^{h}\not\subset H^{1}(\Omega)$.\\
\indent When $n=2$, the Morley-Wang-Xu element space is the Morley element space.\\
\indent The discrete eigenvalue problem reads: Find
$(\lambda_{h},u_{h})\in R\times S^{h}$ with $\|u_{h}\|_{b}=1$ such
that
\begin{eqnarray}\label{2.4}
a_{h}(u_{h},v)=\lambda_{h} b(u_{h},v),~~~\forall v\in S^{h},
\end{eqnarray}
where
\begin{eqnarray*}
&&a_{h}(u_{h},v)= \sum\limits_{\kappa\in\pi_{h}}
\int\limits_{\kappa} ( \sum\limits_{i,j=1}^{n}a_{ij}
\frac{\partial^{2}u_{h}}{\partial x_{i}\partial
x_{j}}\frac{\partial^{2}v}{\partial x_{i}\partial x_{j}}
+\sum\limits_{i,j=1}^{n}b_{ij} \frac{\partial u_{h}}{\partial
x_{i}}\frac{\partial v}{\partial x_{j}}
 +\beta u_{h} v) dx.
\end{eqnarray*}
\indent Let
\begin{eqnarray*}
&&\|v\|_{m,p,h}^{p}=\sum\limits_{\kappa\in
\pi_{h}}\|v\|_{m,p,\kappa}^{p},~~~
|v|_{m,p,h}^{p}=\sum\limits_{\kappa\in
\pi_{h}}|v|_{m,p,\kappa}^{p},\\
&& \|v\|_{m,h}=\|v\|_{m,2,h} ~~~|v|_{m,h}=|v|_{m,2,h},~~~ m=0,1,2.
\end{eqnarray*}
\indent From Lemma 8 in \cite{wm} we know that $|\cdot|_{2,h}$ is
equivalent to $\|\cdot\|_{2,h}$, $\|\cdot\|_{2,h}$ is a norm in
$S^{h}$, and $a_{h}(\cdot,\cdot)$ is a uniformly $S^{h}$-elliptic
bilinear form,
and $\|\cdot\|_{h}=a_{h}(\cdot,\cdot)^{\frac{1}{2}}$ is a norm in $S^{h}$.\\
\indent With regard to the error estimate of the Morley-Wang-Xu
element approximation for biharmonic equations we refer to
\cite{lascaux,mao,shi1,shi2,wm}, and as for the biharmonic
eigenvalue problems we refer to \cite{hu1,rr2} where the work can
be extended to
(\ref{2.1})-(\ref{2.2}).\\
\indent Define
\begin{eqnarray*}
&&P_{F}^{0}f=\frac{1}{meas(F)}\int\limits_{F}fds,~~~R_{F}^{0}f=f-P_{F}^{0}f,\\
&&P_{\kappa}^{0}f=\frac{1}{meas(\kappa)}\int\limits_{\kappa}fdx,~~~R_{\kappa}^{0}f=f-P_{\kappa}^{0}f,
\end{eqnarray*}
where $\kappa\in\pi_{h}$ and $F\in\varepsilon_{h}$.\\
\indent When $\hat{w}\in W_{1,\iota}(\hat{\kappa})$ and $1\leq
g<\frac{(n-1)\iota}{n-\iota}$, by the trace theorem we get
$W_{1,\iota}(\hat{\kappa})\hookrightarrow L_{g}(\partial\hat{\kappa})$.
Thus we can deduce
\begin{eqnarray}\label{2.5}
\int\limits_{\partial\kappa} w^{g}ds\lesssim
h_{\kappa}^{n-\frac{gn}{\iota}-1}\|w\|_{0,\iota,\kappa}^{g}+h_{\kappa}^{g+n-\frac{gn}{\iota}-1}|w|_{1,\iota,\kappa}^{g},~~~\forall
\kappa\in \pi_{h},
\end{eqnarray}
where $\hat{\kappa}$ is a reference element, $\kappa$ and
$\hat{\kappa}$ are affine-equivalent.\\
We define
$$E=(\frac{2n-2}{n}, 2].$$
And we suppose that $p\in E$, $\frac{1}{p'}=1-\frac{1}{p}$.\\
\indent Consider the following associated source problem (\ref{2.6}) and discrete source
problem (\ref{2.7}): Find $w\in H_{0}^{2}(\Omega)$, such that
\begin{eqnarray}\label{2.6}
a(w,v)=b(f,v),~~~\forall v\in H_{0}^{2}(\Omega).
\end{eqnarray}
Find $w_{h}\in S^{h}$,  such that
\begin{eqnarray}\label{2.7}
a_{h}(w_{h},v)=b(f,v),~~~\forall v\in S^{h}.
\end{eqnarray}
Define the consistency term
\begin{eqnarray*}
E_{h}(w,v_{h})=b(f,v_{h})-a_{h}(w,v_{h}),~~~\forall v_{h}\in
S^{h}+H_{0}^{2}(\Omega).
\end{eqnarray*}
\indent Using the proof method in \cite{shi1,wm} we obtain the following error estimate of the consistency term.\\
\indent{\bf Lemma 2.1.}~~Let $w\in W_{3,p}(\Omega)$ be the solution
of (\ref{2.6}), and $p\in E$, then $\forall v_{h}\in
S^{h}+H_{0}^{2}(\Omega)$ there holds
\begin{eqnarray}\label{2.8}
|E_{h}(w,v_{h})|\lesssim
h^{n(\frac{1}{2}-\frac{1}{p})+1}(\|w\|_{3,p}+h^{1-n(\frac{1}{2}-\frac{1}{p})}\|f\|_{b})\|v_{h}\|_{h}.
\end{eqnarray}
\indent{\bf Proof.}~~ For any $v_{h}\in S^{h}$, by Lemma 6 in \cite{wm}(pp. 12,
line 12) we get that there exists a piecewise linear
function $v_{h}^{I}$ on $\pi_{h}$, $v_{h}^{I}\in H_{0}^{1}(\Omega)$,
such that
\begin{eqnarray*}
|v_{h}-v_{h}^{I}|_{m,2,\kappa}\lesssim
h_{\kappa}^{2-m}\sum\limits_{\kappa'\in\pi_{h}(\kappa)}|v_{h}|_{2,2,\kappa'},~~~m=0,1,
\end{eqnarray*}
which together with the inverse inequality yields
\begin{eqnarray*}
&&|v_{h}-v_{h}^{I}|_{m,p',\kappa}\lesssim
h_{\kappa}^{n(\frac{1}{p'}-\frac{1}{2})}
|v_{h}-v_{h}^{I}|_{m,2,\kappa}\\
&&~~~\lesssim
h_{\kappa}^{n(\frac{1}{2}-\frac{1}{p})+2-m}\sum\limits_{\kappa'\in\pi_{h}(\kappa)}|v_{h}|_{2,2,\kappa'},~~~m=0,1,
\end{eqnarray*}
thus, by the Jensen's inequality, we get
\begin{eqnarray}\label{2.9}
|v_{h}-v_{h}^{I}|_{m,p',h}\lesssim
h^{n(\frac{1}{2}-\frac{1}{p})+2-m} |v_{h}|_{2,h},~~~m=0,1.
\end{eqnarray}
For any $v\in H_{0}^{2}(\Omega)$, since $p\in E$, from the
interpolation error estimates we know that there exists a piecewise
linear function $v^{I}\in H_{0}^{1}(\Omega)$ such that (\ref{2.9})
is valid,
and thus, $\forall v_{h}\in S^{h}+H_{0}^{2}(\Omega)$ there exists a $v_{h}^{I}\in H_{0}^{1}(\Omega)$ such that (\ref{2.9}) holds.\\
Write
\begin{eqnarray}\label{2.10}
E_{h}(w,v_{h})=b(f,v_{h}-v_{h}^{I})+b(f, v_{h}^{I})-a_{h}(w,v_{h}).
\end{eqnarray}
From (\ref{2.9}) we have
\begin{eqnarray}\label{2.11}
b(f,v_{h}-v_{h}^{I})\lesssim h^{2}\|f\|_{b}|v_{h}|_{2,h}.
\end{eqnarray}
Since $v_{h}^{I}\in H_{0}^{1}(\Omega)$, by the Green's formula we deduce
\begin{eqnarray}\label{2.12}
&&~~~b(f,v_{h}^{I})-a_{h}(w,v_{h})\nonumber\\
&&=\sum\limits_{\kappa\in\pi_{h}}\int\limits_{\kappa} \{ Lw
v_{h}^{I}- \sum\limits_{i,j=1}^{n}a_{ij}
\frac{\partial^{2}w}{\partial x_{i}\partial
x_{j}}\frac{\partial^{2}v_{h}}{\partial x_{i}\partial x_{j}}
-\sum\limits_{i,j=1}^{n}b_{ij} \frac{\partial w}{\partial
x_{i}}\frac{\partial v_{h}}{\partial x_{j}}
-\beta w v_{h}\}dx\nonumber\\
&&=\sum\limits_{\kappa\in\pi_{h}}\int\limits_{\kappa} \{
\sum\limits_{i,j=1}^{n} \frac{\partial}{\partial x_{i}}
(a_{ij}\frac{\partial^{2}w}{\partial x_{i}\partial x_{j}})
\frac{\partial}{\partial x_{j}}(v_{h}-v_{h}^{I})\nonumber\\
&&~~~-\sum\limits_{i,j=1}^{n}b_{ij} \frac{\partial w}{\partial
x_{i}}\frac{\partial (v_{h}-v_{h}^{I})}{\partial x_{j}}
- \beta w (v_{h}-v_{h}^{I}) \}dx\nonumber\\
&&~~~-\sum\limits_{\kappa\in\pi_{h}} \sum\limits_{1\leq i,j\leq
n}\int\limits_{\partial\kappa} a_{ij}\frac{\partial^{2}w}{\partial
x_{i}\partial x_{j}}\frac{\partial v_{h}}{\partial
x_{j}}\gamma_{i}ds.
\end{eqnarray}
By the H\"{o}lder inequality and (\ref{2.9}) we get
\begin{eqnarray}\label{2.13}
&&~~~|\sum\limits_{\kappa\in\pi_{h}}\int\limits_{\kappa}
 \sum\limits_{i,j=1}^{n}
\frac{\partial}{\partial x_{i}}
(a_{ij}\frac{\partial^{2}w}{\partial x_{i}\partial x_{j}})
\frac{\partial}{\partial
x_{j}}(v_{h}-v_{h}^{I})dx|\nonumber\\
&&\lesssim
\sum\limits_{\kappa\in\pi_{h}}\|w\|_{3,p,\kappa}\|v_{h}-v_{h}^{I}\|_{1,p',\kappa}\lesssim
h^{n(\frac{1}{2}-\frac{1}{p})+1}\|u\|_{3,p}\|v_{h}\|_{h},
\end{eqnarray}
and
\begin{eqnarray}\label{2.14}
|- \sum\limits_{\kappa\in
\pi_{h}}\int\limits_{\kappa}\sum\limits_{i,j=1}^{n}b_{ij}
\frac{\partial w}{\partial x_{i}}\frac{\partial
(v_{h}-v_{h}^{I})}{\partial x_{j}}- \beta w
(v_{h}-v_{h}^{I})dx|\lesssim h \|w\|_{1,2}\|v_{h}\|_{h}.
\end{eqnarray}
It follows from the fact $p\in E$ that
$W_{1,2}(\hat{\kappa})\hookrightarrow L_{p'}(\partial\hat{\kappa})$,
thus from the trace inequality (\ref{2.5}) and the interpolation
error estimate we get
\begin{eqnarray*}
&&\|R_{\kappa}^{0}(a_{ij}\frac{\partial^{2}w}{\partial x_{i}\partial
x_{j}})\|_{0,p,F} \lesssim
h_{\kappa}^{-\frac{1}{p}}\|R_{\kappa}^{0}(a_{ij}\frac{\partial^{2}w}{\partial
x_{i}\partial x_{j}})\|_{0,p,\kappa}
+h_{\kappa}^{1-\frac{1}{p}}\|R_{\kappa}^{0}(a_{ij}\frac{\partial^{2}w}{\partial
x_{i}\partial x_{j}})\|_{1,p,\kappa}\\
&&~~~ \lesssim h_{\kappa}^{1-\frac{1}{p}}\|w\|_{3,p,\kappa},
\end{eqnarray*}
and
\begin{eqnarray*}
\|R_{\kappa}^{0}\frac{\partial v_{h}}{\partial x_{j}}\|_{0,p',F}
\lesssim
h_{\kappa}^{\frac{n}{2}-\frac{n-1}{p}-1}\|R_{\kappa}^{0}\frac{\partial
v_{h}}{\partial x_{j}}\|_{0,2,\kappa}
+h_{\kappa}^{\frac{n}{2}-\frac{n-1}{p}}\|R_{\kappa}^{0}\frac{\partial
v_{h}}{\partial x_{j}}\|_{1,2,\kappa}\lesssim
h_{\kappa}^{\frac{n}{2}-\frac{n-1}{p}}\|v_{h}\|_{2,2,\kappa}.
\end{eqnarray*}
From the above two relations, noting that $\int_{F}[\frac{\partial
v_{h}}{\partial x_{j}}]ds=0$ (see Lemma 4 in \cite{wm}) and
$\|R_{F}^{0}f\|_{0,p,F}\leq 2\|f-v\|_{0,p,F}$ for any $v\in
P_{0}(\kappa)$, we deduce
\begin{eqnarray}\label{2.15}
&&~~~|\sum\limits_{\kappa\in\pi_{h}} \sum\limits_{1\leq i,j\leq
n}\int\limits_{\partial\kappa}
a_{ij}\frac{\partial^{2}w}{\partial x_{i}\partial x_{j}}\frac{\partial v_{h}}{\partial x_{j}}\gamma_{i}ds|\nonumber\\
&&=|\sum\limits_{\kappa\in\pi_{h}} \sum\limits_{1\leq i,j\leq
n}\sum\limits_{F\subset\partial \kappa}\int\limits_{F}
a_{ij}\frac{\partial^{2}w}{\partial x_{i}\partial x_{j}}\frac{\partial v_{h}}{\partial x_{j}}(\gamma_{F})_{i}dF|\nonumber\\
&&=|\sum\limits_{\kappa\in\pi_{h}} \sum\limits_{1\leq i,j\leq
n}\sum\limits_{F\subset\partial \kappa}\int\limits_{F}
R_{F}^{0}(a_{ij}\frac{\partial^{2}w}{\partial x_{i}\partial x_{j}})R_{F}^{0}\frac{\partial v_{h}}{\partial x_{j}}(\gamma_{F})_{i}dF|\nonumber\\
&&\lesssim \sum\limits_{\kappa\in\pi_{h}} \sum\limits_{1\leq i,j\leq
n}\sum\limits_{F\subset\partial \kappa}
|R_{\kappa}^{0}(a_{ij}\frac{\partial^{2}w}{\partial x_{i}\partial
x_{j}})|_{0,p,F}
|R_{\kappa}^{0}\frac{\partial v_{h}}{\partial x_{j}}|_{0,p',F}\nonumber\\
&&\lesssim \sum\limits_{\kappa\in\pi_{h}}
h_{\kappa}^{1-\frac{1}{p}+\frac{n}{2}-\frac{n-1}{p}}\|w\|_{3,p,\kappa}\|v_{h}\|_{2,2,\kappa}\nonumber\\
&&\lesssim h^{n(\frac{1}{2}-\frac{1}{p})+1}\|w\|_{3,p}\|v_{h}\|_{h}.
\end{eqnarray}
Substituting (\ref{2.13}), (\ref{2.14}) and (\ref{2.15}) into
(\ref{2.12}), we get
\begin{eqnarray*}
|b(f, v_{h}^{I})-a_{h}(w,v_{h})|\lesssim
h^{n(\frac{1}{2}-\frac{1}{p})+1}\|w\|_{3,p}\|v_{h}\|_{h}.
\end{eqnarray*}
Substituting (\ref{2.11}) into (\ref{2.10}), and combining the
obtained conclusion with the above relation we get (\ref{2.8}).
~~~$\Box$\\

\indent Define the solution operators $T:L_{2}({\Omega})\to
H_{0}^{2}(\Omega)\subset L_{2}({\Omega})$ and $T_{h}:
L_{2}({\Omega})\to S^{h}$ as follows:
\begin{eqnarray*}
a(Tf,v)&=&b(f,v),~~~ \forall v \in H_{0}^{2}(\Omega),\\
a_{h}(T_{h}f,v)&=&b(f,v),~~~ \forall v \in S^{h}.
\end{eqnarray*}
Then $T, T_{h}: L_{2}(\Omega) \to L_{2}(\Omega)$ are all
self-adjoint
completely continuous operators, and $\|T_{h}-T\|_{b}\to 0~(h\to 0)$.\\
\indent We need the following regularity assumption: $\forall f\in
H^{-1}(\Omega)$, $Tf\in W_{3,q}(\Omega)$, $q\in E$, $Tf\in
H^{2+\sigma}(\Omega)$,  and
\begin{eqnarray}\label{2.16}
\|Tf\|_{3,q}\lesssim \|f\|_{-1},~~~\|Tf\|_{2+\sigma}\leq C
\|f\|_{-1},
\end{eqnarray}
where $\sigma=n(\frac{1}{2}-\frac{1}{q})+1$.\\
\indent From Theorem 2 in \cite{blum}, we get that if $\Omega\subset
R^{2}$ and the inner angle $\omega$ at each singular point  
is smaller than $180^0$, then $q=2$.\\

\indent Let $w\in W_{3,p}(\Omega)$ ($p\in E$) be the solution of
(\ref{2.6}), and let $w_{h}$ be the solution of (\ref{2.7}), then by
the Strang Lemma we have
\begin{eqnarray}\label{2.17}
\|w_{h}-w\|_{h}\lesssim h^{n(\frac{1}{2}-\frac{1}{p})+1}\|w\|_{3,p},
\end{eqnarray}
Further assume that (\ref{2.16}) is valid, then from the
Nitsche-Lascaux-Lesaint Lemma we get
\begin{eqnarray}\label{2.18}
\|w_{h}-w\|_{b}\lesssim
h^{(1-\frac{1}{p}-\frac{1}{q})n+2}\|w\|_{3,p}.
\end{eqnarray}

\indent Using the theory of spectral approximation we get the following lemma.\\
\indent{\bf Lemma 2.2.}~~ Let $(\lambda_{h},u_{h})$ be the jth
eigenpair of (\ref{2.4}) with $\|u_{h}\|_{b}=1$, $\lambda$ be the
jth eigenvalue of (\ref{2.3}), $u$ be the eigenfunction
corresponding to $\lambda$ which is approximated by $u_{h}$, and
$\|u\|_{b}=1$. Suppose that $u\in W_{3,p}(\Omega)$, $p\in E$, and
(\ref{2.16}) holds, then
\begin{eqnarray}\label{2.19}
&&\mid\lambda_{h}-\lambda\mid \lesssim \lambda^{2}
h^{n(1-\frac{1}{p}-\frac{1}{q})+2},\\\label{2.20}
&&\|u_{h}-u\|_{h}\lesssim \lambda
h^{n(\frac{1}{2}-\frac{1}{p})+1},\\\label{2.21}
&&\|u_{h}-u\|_{b}\lesssim \lambda
h^{n(1-\frac{1}{p}-\frac{1}{q})+2}.
\end{eqnarray}
 \indent {\bf Proof.}~~
From the theory of spectral approximation, we have (see, e.g.,
\cite{babuska,rr2}, Lemma 2.3 in \cite{yang1})
\begin{eqnarray}\label{2.22}
&&\|u_{h}-u\|_{b}\lesssim \lambda\|(T-T_{h})u\|_{b},\\
\label{2.23} && \|u_{h}-u\|_{h}=\lambda\|
(T-T_{h})u\|_{h}+R,\\\label{2.24}
 &&|\lambda_{h}-\lambda|\lesssim
\lambda^{2}\|(T-T_{h})u\|_{b},
\end{eqnarray}
where $\mid R\mid\lesssim \|(T-T_{h})u\|_{b}$.\\
Let $f=u$ in (\ref{2.6})-(\ref{2.7}), then $ w=Tu,~ w_h=T_hu$. Thus,
from (\ref{2.17}), (\ref{2.18}) and (\ref{2.23}) we get
(\ref{2.20}), and from (\ref{2.18}) and (\ref{2.22}) we get
(\ref{2.21}).
 Substituting (\ref{2.18}) into (\ref{2.24}), we get (\ref{2.19}).
The proof is completed.~~~$\Box$\\

\indent Noting that $2\geq p\geq q$, $p,q\in E$ and (\ref{2.16}), we
have
\begin{eqnarray}\label{2.25}
&&n(\frac{1}{2}-\frac{1}{p})+1\geq
n(\frac{1}{2}-\frac{1}{q})+1\nonumber\\
&&~~~=\sigma>n(\frac{1}{2}-\frac{n}{2n-2})+1 =\frac{n-2}{2n-2}\geq
0,
\end{eqnarray}
thus,
\begin{eqnarray}\label{2.26}
h^{n(\frac{1}{2}-\frac{1}{p})+1}\leq
h^{n(\frac{1}{2}-\frac{1}{q})+1}=h^{\sigma}.
\end{eqnarray}

\setcounter{section}{2}\setcounter{equation}{0}
\section{Asymptotic lower bounds for eigenvalues}

\indent The identity in the following Lemma 3.1 (see, e.g., Lemma
3.1 in \cite{yang4}, Lemma 3.2 in \cite{yang5}), which is an
equivalent form of the identity in \cite{yang1,zhang} and is a
generalization of the identity in \cite{armentano}, is an important
tool in studying non-conforming element eigenvalue approximations.
\\
\indent{\bf Lemma 3.1.}~~ Let $(\lambda,u)\in \mathbf{R}\times
H_{0}^{2}(\Omega)$ be an eigenpair of (\ref{2.3}),
$(\lambda_{h},u_{h})\in \mathbf{R}\times S^{h}$ be an eigenpair of
(\ref{2.4}), then the following identity is valid:
\begin{eqnarray}\label{3.1}
&&\lambda-\lambda_{h}=\|u-u_{h}\|_{h}^{2}-\lambda_{h}\|u-u_{h}\|_{b}^{2}\nonumber\\
&&~~~~~~-2\lambda_{h} b(u-v, u_{h})+2a_{h}(u-v,u_{h}),~~~\forall~
v\in S^{h}.
\end{eqnarray}

\indent Define an interpolation operator $I_{h}$: First, define
$I_{k}:H^{2}(\kappa)\to P_{2}(\kappa)$ such that
$I_{\kappa}:H^{2}(\kappa)\to P_{2}(\kappa)$,
\begin{eqnarray}\label{3.2}
\int_{l}I_{\kappa}v=\int_{l}v,~~~\int_{F}\frac{\partial
I_{\kappa}v}{\partial\gamma}ds=\int_{F}\frac{\partial
v}{\partial\gamma}ds,
\end{eqnarray}
where $l$ and $F$ denote any vertice and face of $\kappa$
respectively. Next, define
\begin{eqnarray*}
I_{h}v|_{\kappa}=I_{\kappa}v,~~~\forall \kappa\in\pi_{h}.
\end{eqnarray*}

\indent{\bf Lemma 3.2.}~~Let $u\in W_{3,p}(\Omega)$ ($p\in E$), and
$W_{3,p}(\hat{\kappa})\hookrightarrow W_{s,g}(\hat{\kappa})$, then
\begin{eqnarray}\label{3.3}
\|I_{h}u-u\|_{s,g,\kappa}\leq C
h_{\kappa}^{n(\frac{1}{g}-\frac{1}{p})+3-s}|u|_{3,p,\kappa},~~~s=0,1,2.
\end{eqnarray}
where $\hat{\kappa}$ is the reference element, $\kappa$ and
$\hat{\kappa}$ are affine-equivalent.\\
\indent{\bf Proof.}~~The proof is standard, e.g., see \cite{brenner} or
Theorem 15.3 of \cite{ciarlet}.
~~~$\Box$\\

\indent The following weak interpolation estimation plays an crucial role in our analysis.\\
\indent{\bf Lemma 3.3.}~~Let $u\in W_{3,p}(\Omega)$ ($p\in E$), then
\begin{eqnarray}\label{3.4}
&&|a_{h}(u-I_{h}u, u_{h})|\lesssim h
\max\limits_{i,j}(\max\limits_{\kappa}|a_{ij}(x)|_{1,\infty,\kappa})\|u-u_{h}\|_{h}\|u_{h}\|_{h}\nonumber\\
&&~~~~~~+
(h^{n(\frac{1}{2}-\frac{1}{p})+2}|u_{h}-u|_{1,2,h}|u|_{3,p}+h^{2}|u|_{3,p}|u|_{1,p'})
\max\limits_{i,j}\|b_{ij}\|_{0,\infty}.\nonumber\\
&&~~~~~~+\|\beta\|_{0,\infty}h^{3}|u|_{3,p}\|u_{h}\|_{h},\\\label{3.5}
&&b(u-I_{h}u, u_{h}) \lesssim
h^{3}|u|_{3,p}\|\rho(x)\|_{0,\infty}\|u_{h}\|_{h}.
\end{eqnarray}
 \indent{\bf Proof.}~~ From (\ref{3.2}) and the Green's formula we deduce
\begin{eqnarray}\label{3.6}
\int\limits_{F}\frac{\partial(u-I_{h}u)}{\partial
x_{j}}ds=0,~~~\forall F\in \varepsilon.
\end{eqnarray}
Let $I_{0}$ be the piecewise constant projection operator, then by
(\ref{3.6}) we have
\begin{eqnarray*}
&&\sum\limits_{\kappa\in\pi_{h}} \int\limits_{\kappa}
\sum\limits_{i,j=1}^{n}I_{0}a_{ij}
\frac{\partial^{2}(u-I_{h}u)}{\partial x_{i}\partial
x_{j}}\frac{\partial^{2}u_{h}}{\partial x_{i}\partial
x_{j}}dx\nonumber\\
&&~~~ =\sum\limits_{\kappa\in\pi_{h}} \int\limits_{\partial\kappa}
\sum\limits_{i,j=1}^{n}I_{0}a_{ij} \frac{\partial
(u-I_{h}u)}{\partial x_{j}}\frac{\partial^{2}u_{h}}{\partial
x_{i}\partial x_{j}}\gamma_{i}ds =0,
\end{eqnarray*}
thus
\begin{eqnarray}\label{3.7}
&&a_{h}(u-I_{h}u,u_{h})= \sum\limits_{\kappa\in\pi_{h}}
\int\limits_{\kappa} \{ \sum\limits_{i,j=1}^{n}(a_{ij}-I_{0}a_{ij})
\frac{\partial^{2}(u-I_{h}u)}{\partial x_{i}\partial x_{j}}\frac{\partial^{2}u_{h}}{\partial x_{i}\partial x_{j}}\nonumber\\
&&~~~~~~+\sum\limits_{i,j=1}^{n}b_{ij} \frac{\partial (u-I_{h}u)
}{\partial x_{i}}\frac{\partial u_{h}}{\partial x_{j}}+\beta
(u-I_{h}u) u_{h}\} dx.
\end{eqnarray}
Noticing that $|u-I_{h}u|_{2,h}\lesssim \|u-u_{h}\|_{h}$, using the
interpolation error estimates  we deduce
\begin{eqnarray}\label{3.8}
&&|\sum\limits_{\kappa\in\pi_{h}} \int\limits_{\kappa} (
\sum\limits_{i,j=1}^{n}(a_{ij}-I_{0}a_{ij})
\frac{\partial^{2}(u-I_{h}u)}{\partial x_{i}\partial x_{j}}\frac{\partial^{2}u_{h}}{\partial x_{i}\partial x_{j}}dx|\nonumber\\
&&~~~\lesssim h \max\limits_{i,j}(\max\limits_{\kappa}|a_{ij}(x)|_{1,\infty,\kappa})\|u-I_{h}u\|_{h}\|u_{h}\|_{h}\nonumber\\
&&~~~\lesssim h
\max\limits_{i,j}(\max\limits_{\kappa}|a_{ij}(x)|_{1,\infty,\kappa})\|u-u_{h}\|_{h}\|u_{h}\|_{h}.
\end{eqnarray}
Using the H$\ddot{o}$lder inequality and the error estimate of
interpolation, we obtain
\begin{eqnarray}\label{3.9}
&&|\sum\limits_{\kappa}\int\limits_{\kappa}
\sum\limits_{i,j=1}^{n}b_{ij} \frac{\partial ( u-I_{h}u)}{\partial
x_{i}}\frac{\partial u_{h}}{\partial x_{j}} dx|\nonumber\\
&&~~~\lesssim (|u-I_{h}u|_{1,2,h}|u_{h}-u|_{1,2,h}+|u-I_{h}u|_{1,
p,h}|u|_{1,p'})\max\limits_{i,j}\|b_{ij}\|_{0,\infty}
\nonumber\\
&&~~~\lesssim
(h^{n(\frac{1}{2}-\frac{1}{p})+2}|u_{h}-u|_{1,2,h}|u|_{3,p}+h^{2}|u|_{3,p}|u|_{1,p'})
\max\limits_{i,j}\|b_{ij}\|_{0,\infty}.
\end{eqnarray}
From (\ref{2.9}), we get
\begin{eqnarray}\label{3.10}
&&\|u_{h}\|_{0,p'}\leq
\|u_{h}-u_{h}^{I}\|_{0,p'}+\|u_{h}^{I}\|_{0,p'}\lesssim
h^{n(\frac{1}{2}-\frac{1}{p})+2}
|u_{h}|_{2,h}+\|u_{h}^{I}\|_{1,2}\nonumber\\
&&~~~\lesssim h^{n(\frac{1}{2}-\frac{1}{p})+2}
|u_{h}|_{2,h}+\|u_{h}^{I}-u_{h}\|_{1,2,h}+\|u_{h}\|_{1,2,h}\nonumber\\
&&~~~\lesssim \|u_{h}\|_{h}=\lambda_{h}^{\frac{1}{2}}.
\end{eqnarray}
Thus, we obtain
\begin{eqnarray}\label{3.11}
&&\sum\limits_{\kappa\in\pi_{h}} \int\limits_{\kappa}\beta
(u-I_{h}u) u_{h} dx\lesssim
h^{3}\|\beta\|_{0,\infty}|u|_{3,p}\|u_{h}\|_{0,p'}\nonumber\\
&&~~~\lesssim h^{3}\|\beta\|_{0,\infty}|u|_{3,p}\|u_{h}\|_{h}.
\end{eqnarray}
Substituting (\ref{3.8}), (\ref{3.9}) and (\ref{3.11}) into (\ref{3.7}) we get (\ref{3.4}).\\
From the H$\ddot{o}$lder inequality and the error estimate of
interpolation we obtain (\ref{3.5}).
~~~$\Box$\\

\indent The following lemma is another key in our analysis.\\
\indent{\bf Lemma 3.4.}~~Let $w\in W_{3,p}(\Omega)$ ($p\in E$) be
solution of (\ref{2.6}), and let $w_{h}$ be the solution of
(\ref{2.7}), assume that (\ref{2.16}) is valid, then
\begin{eqnarray}\label{3.12}
\|w_{h}-w\|_{b}\lesssim
h^{\sigma}\|w_{h}-w\|_{h}+h^{1+\sigma}\|f\|_{b},
\end{eqnarray}
furthermore, under the conditions of Lemma 2.2, there holds
\begin{eqnarray}\label{3.13}
&&\|u-u_{h}\|_{b}\lesssim h^{\sigma}\|u_{h}-u\|_{h}+h^{1+\sigma}.
\end{eqnarray}
\indent {\bf Proof.}~~ Referring to (\ref{3.7}), we have
\begin{eqnarray*}
&&a_{h}(Tg-I_{h}Tg,w_{h})= \sum\limits_{\kappa\in\pi_{h}}
\int\limits_{\kappa} \{ \sum\limits_{i,j=1}^{n}(a_{ij}-I_{0}a_{ij})
\frac{\partial^{2}(Tg-I_{h}Tg)}{\partial x_{i}\partial x_{j}}\frac{\partial^{2}w_{h}}{\partial x_{i}\partial x_{j}}\nonumber\\
&&~~~~~~+\sum\limits_{i,j=1}^{n}b_{i,j} \frac{\partial (Tg-I_{h}Tg)
}{\partial x_{i}}\frac{\partial w_{h}}{\partial x_{j}}+\beta
(Tg-I_{h}Tg) w_{h}\} dx\\
&&~~~\lesssim h^{1+\sigma}\|g\|_{b}\|w_{h}\|_{h},
\end{eqnarray*}
thus
\begin{eqnarray}\label{3.14}
&&E_{h}(w,Tg-I_{h}Tg)=b(f,Tg-I_{h}Tg)-a_{h}(w-w_{h},Tg-I_{h}Tg)\nonumber\\
&&~~~~~~-a_{h}(w_{h},Tg-I_{h}Tg)\nonumber\\
&&~~~\lesssim
h^{2+\sigma}\|Tg\|_{2+\sigma}\|f\|_{b}+h^{\sigma}\|Tg\|_{2+\sigma}\|w-w_{h}\|_{h}\nonumber\\
&&~~~~~~+h^{1+\sigma}\|g\|_{b}\|w_{h}\|_{h}.
\end{eqnarray}
 By the Nitsche-Lascaux-Lesaint Lemma, (\ref{3.3}), (\ref{3.14}) and (\ref{2.8}), we have
\begin{eqnarray*}
&&\|w-w_{h}\|_{b}\lesssim \|w-w_{h}\|_{h}\sup\limits_{g\in
L_{2}(\Omega), g\not=0}\{\frac{1}{\|g\|_{b}} \|Tg-I_{h}Tg\|_{h}\}
\nonumber\\
&&~~~~~~+\sup\limits_{g\in H, g\not=0}  \{\frac{1}{\|g\|_{b}}
(E_{h}(w,Tg-I_{h}Tg) +E_{h}(Tg,w-w_{h}))\}\nonumber\\
&&~~~\lesssim h^{\sigma}\|w-w_{h}\|_{h}+h^{1+\sigma}\|f\|_{b},
\end{eqnarray*}
i.e., (\ref{3.12}) is valid.\\
\indent From (\ref{3.12}) and (\ref{2.23}) we get
\begin{eqnarray*}
&&\|Tu-T_{h}u\|_{b}\lesssim
h^{\sigma}\|Tu-T_{h}u\|_{h}+h^{1+\sigma}\|u\|_{b}\\
&&~~~ \lesssim
h^{\sigma}\|u_{h}-u\|_{h}+h^{\sigma}\|Tu-T_{h}u\|_{b}+h^{1+\sigma},
\end{eqnarray*}
thus,
\begin{eqnarray*}
\|Tu-T_{h}u\|_{b}\lesssim h^{\sigma}\|u_{h}-u\|_{h}+h^{1+\sigma},
\end{eqnarray*}
which together with (\ref{2.22}) yields (\ref{3.13}).
The proof is completed.~~~$\Box$\\
\indent  The above (\ref{3.13}) is first given in \cite{lin2} for
biharmonic eigenvalue problems (see (60) in
\cite{lin2}) while a detailed proof is not provided. Here, we give a proof of (\ref{3.13}) for more general fourth-order problems (\ref{2.1})-(\ref{2.2}).\\

\indent Based on the above lemmas, we can easily get the following
theorems on the lower bound property of the Morley element
eigenvalues for fourth order elliptic eigenvalue problems in any
dimensions.\\

\indent{\bf Theorem 3.1.}~~ Under the conditions of Lemma 2.2,
suppose that there exists $p_{0}$ satisfying $p<p_{0}<2$,
 $u\in W_{3,p}(\Omega)$, $u\not\in W_{3,p_{0}}(\Omega)$
. And suppose that $\|u_{h}-u\|_{h}\gtrsim h^{1-\delta}$ with
$\delta\in (0,n/p_{0}-n/2)$ be an arbitrarily small constant, then
when $h$ is small enough there holds
\begin{eqnarray}\label{3.15}
\lambda_{h}\leq \lambda.
\end{eqnarray}
\indent{\bf Proof.}~~ Taking $v_{h}=I_{h}u$ in (\ref{3.1}) we get
\begin{eqnarray}\label{3.16}
&&\lambda
-\lambda_{h}=\|u-u_{h}\|_{h}^{2}-\lambda_{h}\|u-u_{h}\|_{b}^{2}-2\lambda_{h}b(u-I_{h}u,
u_{h})\nonumber\\
&&~~~~~~+2a_{h}(u-I_{h}u, u_{h}).
\end{eqnarray}
Next we shall compare the four terms on the right-hand side of (\ref{3.16}).\\
\indent From (\ref{3.13}), we get
\begin{eqnarray}\label{3.17}
&&\|u-u_{h}\|_{b}\lesssim
h^{\sigma}\|u_{h}-u\|_{h}+h^{1+\sigma}\lesssim
h^{\sigma}\|u_{h}-u\|_{h},
\end{eqnarray}
which indicates that in (\ref{3.16}) the second term is a quantity of higher order than the first term.\\
From (\ref{3.5}), we have
\begin{eqnarray*}
&&b(u-I_{h}u, u_{h})\lesssim
h^{3}\|\rho(x)\|_{0,\infty}|u|_{3,p}\|u_{h}\|_{h},
\end{eqnarray*}
which implies that the third term is also a quantity of higher order than the first term.\\
From (\ref{3.4}) and (\ref{2.26}) we have
\begin{eqnarray*}
&&|a_{h}(u-I_{h}u, u_{h})|\lesssim h
\max\limits_{i,j}|a_{ij}(x)|_{1,\infty}\|u-u_{h}\|_{h}\|u_{h}\|_{h}\nonumber\\
&&~~~~~~+ h^{1+\sigma}\|u_{h}-u\|_{h}+
 h^{2}+h^{3}.
\end{eqnarray*}
Thus, the fourth term is quantity of higher order than the first one.\\
Hence, (\ref{3.15}) is valid.
~~~$\Box$\\

\indent{\bf Theorem 3.2.}~~Under the conditions of Lemma 2.2, let
$a_{ij}(x)$ be piecewise constants, $b_{ij}=0$ ($i,j=1,2,\cdots,n$).
And suppose that $\|u_{h}-u\|_{h}\gtrsim h$ and $h$ is small enough,
then there holds
\begin{eqnarray}\label{3.18}
\lambda_{h}\leq \lambda.
\end{eqnarray}
\indent{\bf Proof.}~~We shall compare the four terms on the
right-hand side of (\ref{3.16}).\\
From (\ref{3.13}) and (\ref{3.5}) we know that
the second and the third term are quantities of higher order than the first one, respectively.\\
Noting that $a_{ij}(x)$ is a constant and $b_{ij}=0$, from
(\ref{3.4}) we have
\begin{eqnarray*}
|a_{h}(u-I_{h}u, u_{h})|\lesssim
h^{3}\|\beta\|_{0,\infty}|u|_{3,p}\|u_{h}\|_{h}\lesssim
h^{3}|u|_{3,p}.
\end{eqnarray*}
which indicates that the fourth term is quantity of higher order than the first one.\\
Hence, (\ref{3.18}) is valid.
~~~$\Box$\\

\indent From the proof of Theorem 3.2 we know that the condition, $\|u_{h}-u\|_{h}\gtrsim h$ , can be weakened to $\|u_{h}-u\|_{h}\gtrsim h^{1+\frac{\sigma}{3}}$.\\
\indent{\bf Corollary 3.1.}~~ Under the conditions of Theorem 3.2,
assume that $\|u_{h}-u\|_{h}\gtrsim h^{1+\frac{\sigma}{3}}$ and $h$
is small enough, then there holds
\begin{eqnarray}\label{3.19}
\lambda_{h}\leq \lambda.
\end{eqnarray}

\indent One noticed early that the error of finite element
eigenvalues is relevant to the value of eigenvalues, which means
that the computation will be more difficult when the eigenvalue becomes larger
(see Section 6.3 in \cite{strang}). However, we find when the value of eigenvalues is large,
the lower bound for such eigenvalues can
be obtained using the Morley element.\\

\indent{\bf Theorem 3.3.}~~Under the conditions of Lemma 2.2,
suppose $p=2$ and $\|u-u_{h}\|_{h}\gtrsim h|u|_{3,2}$. Then, for
the eigenvalue $\lambda$ that the value is large, when $h$ is small enough it is valid
that
\begin{eqnarray}\label{3.20}
\lambda_{h}\leq \lambda.
\end{eqnarray}
\indent{\bf Proof.}~~We shall compare the four terms on the
right-hand side of (\ref{3.16}).\\
\indent From (\ref{3.13}) and (\ref{3.5}) we know that the second
and third
term are all quantities of higher order than the first one.\\
A simple calculation shows that
\begin{eqnarray}\label{3.21}
|u|_{1,2}^{2}=\int\limits_{\Omega}\nabla u\cdot\nabla
udx=-\int\limits_{\Omega}u \triangle udx\lesssim
\|u\|_{b}\|u\|_{a}=\lambda^{\frac{1}{2}},
\end{eqnarray}
and
\begin{eqnarray}\label{3.22}
\lambda=\|u\|_{a}^{2}\thickapprox |u|_{2,2}^{2}=\|\triangle
u\|_{0,2}^{2} \lesssim |u|_{1,2}|u|_{3,2}\lesssim
\lambda^{\frac{1}{4}}|u|_{3,2}.
\end{eqnarray}
Substituting (\ref{3.21}) and
$\|u_{h}\|_{h}\thickapprox \lambda^{\frac{1}{2}}$ into (\ref{3.4}),
 we get
\begin{eqnarray*}
&&|a_{h}(u-I_{h}u, u_{h})|\lesssim h
\max\limits_{i,j}(\max\limits_{\kappa}|a_{ij}(x)|_{1,\infty,\kappa})\|u-u_{h}\|_{h}\lambda^{\frac{1}{2}}\nonumber\\
&&~~~~~~+
(h^{2}\|u_{h}-u\|_{h}|u|_{3,2}+h^{2}|u|_{3,2}\lambda^{\frac{1}{4}})
\max\limits_{i,j}\|b_{ij}\|_{0,\infty}
+h^{3}|u|_{3,2}\|\beta\|_{0,\infty}\lambda^{\frac{1}{2}},
\end{eqnarray*}
by the stability assumption and (\ref{3.22}), we have
\begin{eqnarray*}
\|u-u_{h}\|_{h}\gtrsim h|u|_{3,2}\gtrsim h\lambda^{\frac{3}{4}},
\end{eqnarray*}
which indicates that,  for the eigenvalue $\lambda$ that the value is large, when $h$
is small enough the absolute value of the fourth term is smaller
than that of
the first one.\\
Thus, (\ref{3.20}) is valid.~~~$\Box$\\

\indent {\bf Remarks 3.1.}~~From \cite{babuska2,hu1,lin3,widlund} we
know that the saturation condition $\|u-u_{h}\|_{h}\gtrsim
h^{t_{0}}$ holds on the quasi-uniform mesh $\pi_{h}$, where $t_{0}$
is the singularity exponent of the eigenfunction $u$. However, this
condition isn't valid on adaptive meshes with local refinement.
Inspired by \cite{lin2,yang2} we change the condition
$\|u-u_{h}\|_{h}\geq Ch^{t_{0}}$ into $\|u_{h}-u\|_{h}\gtrsim
h^{1-\delta}$ with $\delta\in (0,n/p_{0}-n/2)$ be an arbitrary small
constant in Theorem 3.1, and into $\|u-u_{h}\|_{h}\gtrsim h$ in
Theorem 3.2, respectively. The modified conditions are valid not
only for the quasi-uniform mesh but also for many kinds of adaptive
meshes.
From (8.11) in \cite{hu1} we know that the saturation
condition
$\|u-u_{h}\|_{h}\gtrsim h|u|_{3,2}$ holds in Theorem 3.3.

\indent {\bf Remarks 3.2.}~~When the conditions of Theorem 3.1 or
Theorem 3.2 hold, from Theorem 3.1 or Theorem 3.2 we know that
$\|u-u_{h}\|_{h}^{2}$ is the dominant term among the four terms on
the right-hand side of (\ref{3.16}), i.e.,
\begin{eqnarray}\label{3.23}
\lambda -\lambda_{h}\backsimeq\|u-u_{h}\|_{h}^{2},
\end{eqnarray}
which together with (\ref{2.20}) we obtain
\begin{eqnarray}\label{3.24}
0< \lambda-\lambda_{h} \lesssim h^{(1-\frac{2}{p})n+2}.
\end{eqnarray}
When $p>q$, (\ref{3.24}) is slightly better than (\ref{2.19}).

\indent {\bf Remarks 3.3.}~~Referring to \cite{hu1,luo,yang4} we can
also use conforming finite elements to do the postprocessing to get
the upper bound of eigenvalues.

\indent {\bf Remarks 3.4.}~~We see that $|a_{h}(u-I_{h}u, u_{h})|$
plays a crucial role in Theorem 3.1-Theorem 3.3, and we can get the
following result: Under the conditions of Lemma 2.2, let
$\varepsilon_{0}$ be a upper bound of $|a_{h}(u-I_{h}u, u_{h})|$,
and suppose that $\|u_{h}-u\|_{h}\gtrsim h$ and $h$ is small enough,
then there holds
\begin{eqnarray}\label{3.25}
\lambda_{h}-\varepsilon_{0}\leq \lambda.
\end{eqnarray}

\setcounter{section}{3}\setcounter{equation}{0}
\section{A posteriori error estimates for eigenvalues }

\indent \cite{dai1,yang5} gave the relationship between the
conforming/nonconforming finite element eigenvalue approximation and
the associated conforming/nonconforming finite element boundary
value approximation (see Theorem 3.1 and Lemma 5.1 in \cite{dai1},
Theorem 3.1 in \cite{yang5}). The following Lemma 4.1 is given in
\cite{yang5} (see Theorem 3.1 in \cite{yang5}).\\

\indent Consider the source problem (\ref{2.6}) associated with
(\ref{2.3}) with $f=\lambda_{h}u_{h}$, whose generalized solution is
$w=\lambda_{h}Tu_{h}$ and the Morley element approximation is
$w_{h}=\lambda_{h}T_{h}u_{h}=u_{h}$.\\
\indent{\bf Lemma 4.1.}~~ Let $(\lambda_{h},u_{h})$ be the jth
eigenpair of (\ref{2.4}) with $\|u_{h}\|_{b}=1$, $\lambda$ be the
jth eigenvalue of (\ref{2.3}), then there exists an eigenfunction
$u$  corresponding to $\lambda$ with $\|u\|_{b}=1$, such that
\begin{eqnarray}\label{s4.1}
\|u_{h}-u\|_{h}&=&\lambda_{h}\|Tu_{h}-T_{h}u_{h}\|_{h}+R,
\end{eqnarray}
where $\mid R\mid \lesssim\|(T-T_{h})u\|_{b}$.\\
\indent {\bf Proof.}~~From the definition of $T$ and $T_{h}$ we
derive
\begin{eqnarray}\label{s4.2}
u-u_{h}&=&\lambda Tu-\lambda_{h}T_{h}u_{h}\nonumber\\
&=&\lambda Tu-\lambda_{h}Tu+\lambda_{h}Tu-\lambda_{h}Tu_{h}+\lambda_{h}Tu_{h}-\lambda_{h}T_{h}u_{h}\nonumber\\
&=&(\lambda
-\lambda_{h})Tu+\lambda_{h}T(u-u_{h})+\lambda_{h}(T-T_{h})u_{h}.
\end{eqnarray}
Denote
\begin{eqnarray*}
R=\|u-u_{h}\|_{h}-\lambda_{h} \|(T-T_{h})u_{h}\|_{h}.
\end{eqnarray*}
Using the triangle inequality, (\ref{s4.2}), (\ref{2.22}) and
(\ref{2.24}) we deduce
\begin{eqnarray*}
&&|R|=|\|u-u_{h}\|_{h}-\lambda_{h} \|(T-T_{h})u_{h}\|_{h}|\nonumber\\
&&~~~\leq\|u-u_{h}-\lambda_{h}(T-T_{h})u_{h}\|_{h}\\
&&~~~=\|(\lambda
-\lambda_{h})Tu+\lambda_{h}T(u-u_{h})\|_{h}\nonumber\\
&&~~~\lesssim |\lambda_{h}-\lambda|+
\|u-u_{h}\|_{b}\lesssim\|(T-T_{h})u\|_{b},
\end{eqnarray*}
which proves (\ref{s4.1}).
~~~$\Box$\\

\indent \cite{yang5} also pointed out that one can use Lemma 3.1 to obtain the a posteriori error estimator
of nonconforming finite element eigenvalues (see Lemma 3.2 and Section 4.2 in \cite{yang5}). Hence we have the following Theorem 4.1:\\

\indent{\bf Theorem 4.1.}~~Under the conditions of Theorem 3.1 or
Theorem 3.2, it is valid that
\begin{eqnarray}\label{4.3}
&&\|u_{h}-u\|_{h}\backsimeq
\|T(\lambda_{h}u_{h})-T_{h}(\lambda_{h}u_{h})\|_{h},\\\label{4.4}
&&\lambda-\lambda_{h}\backsimeq
\|T(\lambda_{h}u_{h})-T_{h}(\lambda_{h}u_{h})\|_{h}^{2}.
\end{eqnarray}
\indent {\bf Proof.}~~ From (\ref{2.23}) and (\ref{3.12}) we know
that
 under the conditions of Theorem 3.1 or Theorem 3.2,
$\|Tu-T_{h}u\|_{b}$ is a quantity of higher order than
$\|u_{h}-u\|_{h}$. So, from (\ref{s4.1}) we get (\ref{4.3}). From
Theorem 3.1 or Theorem 3.2 we know that $\|u-u_{h}\|_{h}^{2}$ is the
dominant term in (\ref{3.16}), i.e.,
\begin{eqnarray}\label{4.5}
\lambda -\lambda_{h}\backsimeq\|u-u_{h}\|_{h}^{2},
\end{eqnarray}
which together with (\ref{4.3}) yields (\ref{4.4}).
~~~$\Box$\\
\indent Theorem 4.1 tells us that the error estimates of the Morley
element eigenvalue and eigenfunction are reduced to the error
estimates of the Morley element solution $w_{h}=u_{h}$ of the
associated source problem with the right-hand side
$f=\lambda_{h}u_{h}$. Thus, the a posteriori error estimator of the
Morley element solution for source problem
becomes the a posteriori error estimator of the Morley element eigenfunction and eigenvalue.\\
\indent The a posteriori error estimates for the Morley plate
bending element have been studied, for example, see \cite{bl,hu2}.
In the following we introduce the work in \cite{hu2}:\\
\indent Consider the source problem (\ref{2.6}) and discrete source
problem (\ref{2.7}) where $n=2$, $a_{ij}(x)$ is a constant,
$b_{ij}=0$, and $\beta=0$.\\
\indent Given any $F\in \varepsilon_{h}$ with the length $h_{F}=|F|$,
let $\gamma_{F}=(\gamma_{1}, \gamma_{2})$ be a fixed unit normal and
$\nu_{F}=(-\gamma_{2}, \gamma_{1})$ be the tangential vector.\\
\indent Hu and Shi \cite{hu2} defined the following estimator:
\begin{eqnarray}\label{4.6}
&&\eta_{h}(f,w_{h},\kappa)^{2}=h_{\kappa}^{4}\|f\|_{0,2,\kappa}^{2}\nonumber\\
&&~~~~~~+\sum\limits_{F\in
\varepsilon_{h}\cap\partial\kappa}h_{F}\|[\frac{1}{2}\{\nabla(\nabla
w_{h}) +\nabla(\nabla
w_{h})^{T}\}\nu_{F}]\|_{0,2,F}^{2},\\\label{4.7}
&&\eta_{h}(f,w_{h},\Omega)^{2}=\sum\limits_{\kappa\in\pi_{h}}\eta_{h}(f,w_{h},\kappa)^{2},
\end{eqnarray}
and proved the following lemma.\\
\indent{\bf Lemma 4.2.}~~Let $w$ be the solution of (\ref{2.6}), and
$w_{h}$ be the solution of (\ref{2.7}). Then
\begin{eqnarray}\label{4.8}
\|w_{h}-w\|_{h}\thickapprox\eta_{h}(f,w_{h},\Omega).
\end{eqnarray}

\indent Combining Theorem 4.1 and Lemma 4.2 we get:\\
\indent{\bf Theorem 4.2.}~~Suppose that $n=2$, $a_{ij}(x)$
 are constants, $b_{ij}=0$ and $\beta=0$. Then under the conditions
of Theorem 3.2 there holds
\begin{eqnarray}\label{4.11}
&&\|u_{h}-u\|_{h}\thickapprox
\eta_{h}(\lambda_{h}u_{h},u_{h},\Omega),\\\label{4.12}
&&|\lambda_{h}-\lambda|\thickapprox
\eta_{h}(\lambda_{h}u_{h},u_{h},\Omega)^{2}.
\end{eqnarray}

\indent From the above (\ref{4.8}) and the proof of Theorem 3.2 in
\cite{yang2}, Shen \cite{shen} also gave (\ref{4.11})-(\ref{4.12}).
The feature of our work is to point out that (\ref{4.11})-(\ref{4.12})
are valid under the conditions of Theorem 3.2, i.e., on adaptive meshes.\\

 \indent Lemma 4.2 and Theorem 4.2 can be extended
to the case of $a_{ij}(x)$ and $b_{ij}=\tau$ are constants, and
$\beta=0$. In this case $\eta_{h}(f,w_{h},\kappa)$ and
$\eta_{h}(f,w_{h},\Omega)$ need to be modified as follows.\\

\begin{eqnarray}\label{4.13R}
&&\eta_{h}(f,w_{h},\kappa)^{2}=h_{\kappa}^{4}\|f\|_{0,2,\kappa}^{2}+h_{\kappa}^{2}\tau |w_{h}|_{1,\kappa}^{2}\nonumber\\
&&~~~~~~+\sum\limits_{F\in
\varepsilon_{h}\cap\partial\kappa}h_{F}\|[\frac{1}{2}\{\nabla(\nabla
w_{h}) +\nabla(\nabla
w_{h})^{T}\}\nu_{F}]\|_{0,2,F}^{2},\\\label{4.14R}
&&\eta_{h}(f,w_{h},\Omega)^{2}=\sum\limits_{\kappa\in\pi_{h}}\eta_{h}(f,w_{h},\kappa)^{2}.
\end{eqnarray}

\indent Using the a posteriori error estimates and consulting the
existing standard algorithm (see, e.g., Algorithm C in \cite{dai1}),
we obtain the following adaptive algorithm of the Morley element for the vibration problem of a clamped plate:\\
\indent{\bf Algorithm 1}\\
Choose parameter $0<\theta<1$.\\
\indent{\bf Step 1.}~Pick any initial mesh $\pi_{h_{0}}$ with mesh size $h_{0}$.\\
\indent{\bf Step 2.}~Solve (\ref{2.4}) on $\pi_{h_{0}}$ for discrete solution $(\lambda_{h_{0}}, u_{h_{0}})$.\\
\indent{\bf Step 3.}~Let $l=0$.\\
\indent{\bf Step 4.}~Compute the local indicators ${\eta}_{h_{l}}^{2}(\lambda_{h_{l}}u_{h_{l}}, u_{h_{l}},\kappa)$.\\
\indent{\bf Step 5.}~Construct
$\widehat{\pi}_{h_{l}}\subset\pi_{h_{l}}$ by {\bf Marking Strategy
E}
and parameter $\theta$.\\
\indent{\bf Step 6.}~Refine $\pi_{h_{l}}$ to get a new mesh
$\pi_{h_{l+1}}$
by Procedure ${\bf Refine}$.\\
\indent{\bf Step 7.}~Solve (\ref{2.4}) on $\pi_{h_{l+1}}$ for discrete solution $(\lambda_{h_{l+1}}, u_{h_{l+1}})$.\\
\indent{\bf Step 8.}~Let $l=l+1$ and go to Step 4.\\

\indent {\bf Marking Strategy E}\\
\indent Given parameter $0<\theta<1$:\\
\indent{\bf Step 1.}~~Construct a minimal subset
$\widehat{\pi}_{h_{l}}$ of $\pi_{h_{l}}$ by selecting some elements
in $\pi_{h_{l}}$ such that
\begin{eqnarray*}
\sum\limits_{\kappa\in
\widehat{\pi}_{h_{l}}}{\eta}_{h_{l}}^{2}(\lambda_{h_{l}}u_{h_{l}},
u_{h_{l}}, \kappa) \geq
\theta{\eta}_{h_{l}}^{2}(\lambda_{h_{l}}u_{h_{l}},
u_{h_{l}},\Omega).
\end{eqnarray*}
\indent{\bf Step 2.}~~Mark all the elements in
$\widehat{\pi}_{h_{l}}$.

\setcounter{section}{4}\setcounter{equation}{0}
\section{Numerical experiments}

\indent Consider the vibration problem of a clamped plate under
tension:
\begin{eqnarray}\label{5.1}
&&\Delta^{2}u-\tau\Delta u=\lambda u,~~~in\Omega,\\\label{5.2}
 &&\frac{\partial u}{\partial
\gamma}=0,~~~u=0,~~~on\partial\Omega,
\end{eqnarray}
where $\Omega\subset R^{2}$ and $\tau$ is the tension  coefficient.
When $\tau=0$, (\ref{5.1})-(\ref{5.2}) is the vibration problem of a clamped plate without tension.\\
\indent The weak form of (\ref{5.1})-(\ref{5.2}) and the Morley
element approximation is
(\ref{2.3}) and (\ref{2.4}), respectively, with $n=2$, $a_{ij}=1$, $b_{ij}=\tau$, $\beta=0$ and $\rho=1$.\\
\indent Weinstein and Chien \cite{weinstein} once obtained lower
bounds of eigenvalues for this problem using relaxed boundary
conditions in 1943. In section 3 and 4, we analyze the asymptotic
lower bounds property and a posteriori error of the Morley element
eigenvalues for this problem. Here we provide the numerical results.
We compute the first two eigenvalues $\lambda_{j}$, $j=1,2$.\\
\indent We show the subdivision way for generating initial
triangulation for
the given unit square and L-shaped domain in Fig. 1.\\
\indent We use the Morley element on uniform triangle meshes and
the conforming Bogner-Fox-Schmit rectangle element (BFS element) on
a uniform rectangle meshes to compute, respectively.
The numerical results $\lambda_{j,h}^{M}$ and $\lambda_{j,h}^{BFS}$ (j=1,2) are listed in Tables 1-6.\\
\indent When $\Omega=(0,1)\times (0,\frac{1}{2})\cup
(0,\frac{1}{2})\times (\frac{1}{2}, 1)$ is an L-shaped domain, we
use Algorithm $1$ with the Morley element to compute
$\lambda_{h_{l}}$ on adaptive meshes with the initial mesh diameter
$h=\frac{\sqrt{2}}{32}$. We take $\theta=0.25$ and our program is
compiled using the package of Chen \cite{chenl}.
The numerical results $\lambda_{j,h_{l}}^{M}$ (j=1,2) are shown in Tables 7-9.\\
\indent In Tables 7-9, we use Ndof to denote the number of degrees of freedom.\\
\indent In Tables 1-6, the discrete eigenvalues
$\lambda_{j,h}^{BFS}$ are computed as upper bound.\\
\indent When $\tau=0$, in Tables 1, 4 and 7, we see that the
discrete eigenvalues $\lambda_{j,h}^{M}$ monotonically increase
stably, and the Morley
elements has already given the lower bounds of eigenvalues, which coincide with Theorem 3.2.\\
\indent When $\tau>0$, the eigenvalues of this problem are very
large. In Tables 5, 6, 8 and 9, we see that the discrete eigenvalues
$\lambda_{j,h}^{M}$ monotonically increase stably, and the Morley
elements has already given the lower bounds of eigenvalues, which
coincide with Theorem 3.1. And the numerical results in Tables 2 and
3 coincide with Theorem 3.3.
\\
\indent When $\Omega$ is an L-shaped domain, we depict the adaptive
meshes at the $l$th iteration in Figs.2-4, and the a posteriori
error
indicator curve in Fig.5.\\
\indent Comparing the results in Tables 4-6 and 7-9, we see that the
approximations on adaptive meshes are more precise than those on
uniform meshes.\\
\indent From Fig.5 it can be seen that the a posteriori error
indicator curves are nearly parallel to a line with slope $-1$,
which shows that the adaptive algorithm based on the a posteriori
error estimators (\ref{4.6})-(\ref{4.7}) and (\ref{4.13R})-(\ref{4.14R}) achieves about the convergence rate of $\mathcal{O}(\frac{1}{Ndof})$ and is successful.\\

\begin{figure}
  \centering
  \includegraphics[width=4cm]{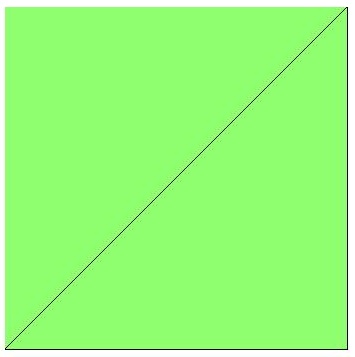}
  \includegraphics[width=4cm]{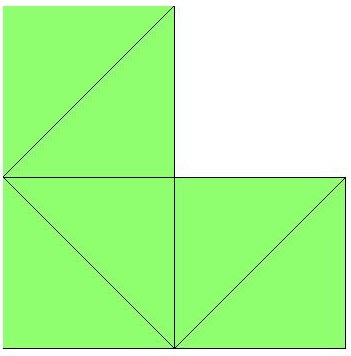}\\
  \caption{Initial meshs}
\end{figure}

\begin{figure}
  \centering
  \includegraphics[width=4cm]{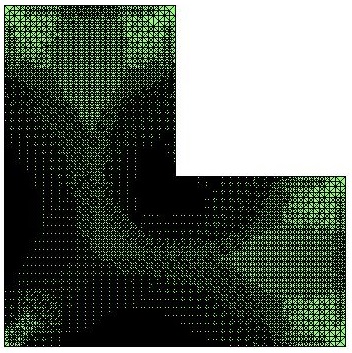}
  \includegraphics[width=4cm]{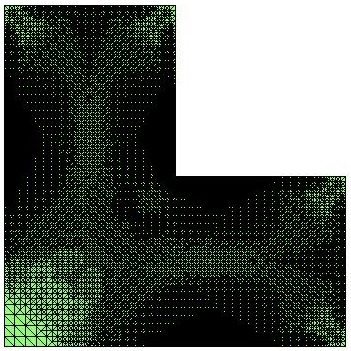}\\
  \caption{$\tau=0$, adaptive meshes of the 27th iteration for the 1st(left) eigenvalue,
  and adaptive meshes of the 24th iteration for the 2nd(right) eigenvalue,
 on the L-shaped domain.}
\end{figure}
\begin{figure}
  \centering
  \includegraphics[width=4cm]{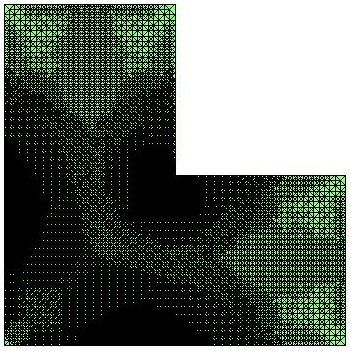}
  \includegraphics[width=4cm]{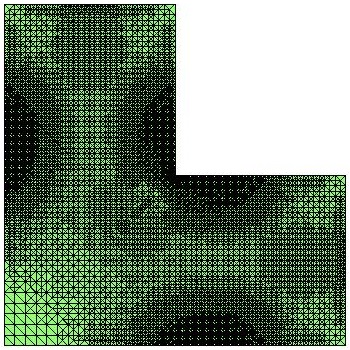}\\
  \caption{$\tau=10$,
  adaptive meshes of the 27th iteration for the 1st(left) eigenvalue,
  and adaptive meshes of the 24th iteration for the 2nd(right) eigenvalue,
 on the L-shaped domain.}
\end{figure}
\begin{figure}
  \centering
  \includegraphics[width=4cm]{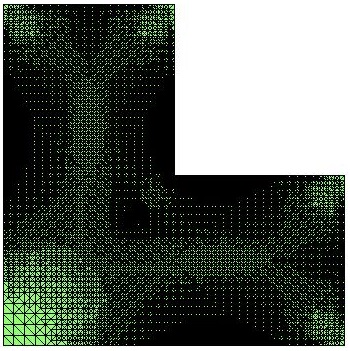}
  \includegraphics[width=4cm]{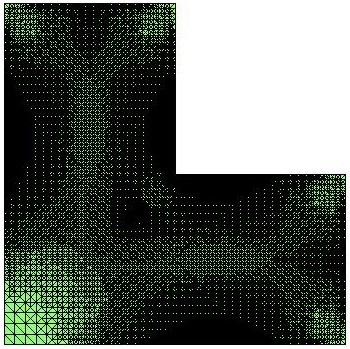}\\
  \caption{$\tau=100$,  adaptive meshes of the 26th iteration for the 1st(left) eigenvalue,
  and adaptive meshes of the 24th iteration for the 2nd(right) eigenvalue,
 on the L-shaped domain.}
\end{figure}

\begin{table}
\small \centering {\bf Table 1: $\tau=0$, the 1st and the 2nd
eigenvalue on the unit square}
\begin{tabular}{ccccccc}\hline
h&Ndof&$\lambda_{1,h}^{M}$&$\lambda_{2,h}^{M}$&Ndof&$\lambda_{1,h}^{BSF}$&$\lambda_{2,h}^{BSF}$\\\hline%
$\frac{\sqrt{2}}{4}$ & 49 & 691.358  & 2068.884  & 36 & 1300.126  & 5480.858 \\
$\frac{\sqrt{2}}{8}$ & 225 & 1049.963  & 3777.006  & 196 & 1295.340  & 5393.253 \\
$\frac{\sqrt{2}}{16}$ & 961 & 1221.316  & 4850.316  & 900 & 1294.963  & 5387.100 \\
$\frac{\sqrt{2}}{32}$ & 3969 & 1275.511  & 5239.489  & 3844 & 1294.936  & 5386.685 \\
$\frac{\sqrt{2}}{64}$ & 16129 & 1290.009  & 5348.921  & 15876 & 1294.934  & 5386.658 \\
$\frac{\sqrt{2}}{128}$ & 65025 & 1293.698  & 5377.161  & 64516 & 1294.934  & 5386.657 \\
$\frac{\sqrt{2}}{256}$ & 261121 & 1294.625  & 5384.279  & ---- & ---- & ----\\
  \hline
\end{tabular}
\end{table}

\begin{table}
\small \centering {\bf Table 2: $\tau=10$, the 1st and the 2nd
eigenvalue on the unit square}
\begin{tabular}{ccccccc}\hline
h&Ndof&$\lambda_{1,h}^{M}$&$\lambda_{2,h}^{M}$&Ndof&$\lambda_{1,h}^{BSF}$&$\lambda_{2,h}^{BSF}$\\\hline%
$\frac{\sqrt{2}}{4}$ & 49 & 817.866  & 2340.295  & 36 & 1540.868  & 6059.160 \\
$\frac{\sqrt{2}}{8}$ & 225 & 1232.645  & 4166.902  & 196 & 1534.567  & 5964.155 \\
$\frac{\sqrt{2}}{16}$ & 961 & 1441.896  & 5352.314  & 900 & 1534.071  & 5957.503 \\
$\frac{\sqrt{2}}{32}$ & 3969 & 1509.595  & 5790.254  & 3844 & 1534.036  & 5957.053 \\
$\frac{\sqrt{2}}{64}$ & 16129 & 1527.828  & 5914.201  & 15876 & 1534.033  & 5957.024 \\
$\frac{\sqrt{2}}{128}$ & 65025 & 1532.476  & 5946.243  & 64516 & 1534.033  & 5957.022 \\
$\frac{\sqrt{2}}{256}$ & 261121 & 1533.643  & 5954.323  & ---- & ---- & ----\\
  \hline
\end{tabular}
\end{table}

\begin{table}
\small \centering {\bf Table 3: $\tau=100$, the 1st and the 2nd
eigenvalue on the unit square}
\begin{tabular}{ccccccc}\hline
h&Ndof&$\lambda_{1,h}^{M}$&$\lambda_{2,h}^{M}$&Ndof&$\lambda_{1,h}^{BSF}$&$\lambda_{2,h}^{BSF}$\\\hline%
$\frac{\sqrt{2}}{4}$ & 49 & 1842.476  & 4567.218  & 36 & 3660.464  & 11228.383 \\
$\frac{\sqrt{2}}{8}$ & 225 & 2684.564  & 7340.279  & 196 & 3624.418  & 11030.745 \\
$\frac{\sqrt{2}}{16}$ & 961 & 3290.809  & 9629.657  & 900 & 3620.928  & 11014.525 \\
$\frac{\sqrt{2}}{32}$ & 3969 & 3528.551  & 10614.330  & 3844 & 3620.671  & 11013.347 \\
$\frac{\sqrt{2}}{64}$ & 16129 & 3596.929  & 10909.494  & 15876 & 3620.654  & 11013.269 \\
$\frac{\sqrt{2}}{128}$ & 65025 & 3614.676  & 10987.054  & 64516 & 3620.653  & 11013.264 \\
$\frac{\sqrt{2}}{256}$ & 261121 & 3619.156  & 11006.694  & ---- & ---- & ----\\
  \hline
\end{tabular}
\end{table}

\begin{table}
\small \centering {\bf Table 4: $\tau=0$, the 1st and the 2nd
eigenvalue on the L-shaped domain}
\begin{tabular}{ccccccc}\hline
h&Ndof&$\lambda_{1,h}^{M}$&$\lambda_{2,h}^{M}$&Ndof&$\lambda_{1,h}^{BSF}$&$\lambda_{2,h}^{BSF}$\\\hline%
$\frac{\sqrt{2}}{4}$ & 33 & 2026.507  & 3077.627  & 20 & 7571.752  & 11513.975 \\
$\frac{\sqrt{2}}{8}$ & 161 & 3897.606  & 6579.447  & 132 & 6999.898  & 11107.297 \\
$\frac{\sqrt{2}}{16}$ & 705 & 5443.844  & 9332.592  & 644 & 6835.442  & 11062.761 \\
$\frac{\sqrt{2}}{32}$ & 2945 & 6223.547  & 10541.203  & 2820 & 6765.112  & 11056.385\\
$\frac{\sqrt{2}}{64}$ & 12033 & 6523.545  & 10916.025  & 11780 & 6732.515  & 11055.009 \\
$\frac{\sqrt{2}}{128}$ & 48641 & 6632.571  & 11018.170  & 48132 & 6717.205  & 11054.646 \\
$\frac{\sqrt{2}}{256}$ & 195585 & 6673.866  & 11045.020  & ---- & ---- & ----\\
  \hline
\end{tabular}
\end{table}

\begin{table}
\small \centering {\bf Table 5: $\tau=10$, the 1st and the 2nd
eigenvalue on the L-shaped domain}
\begin{tabular}{ccccccc}\hline
h&Ndof&$\lambda_{1,h}^{M}$&$\lambda_{2,h}^{M}$&Ndof&$\lambda_{1,h}^{BSF}$&$\lambda_{2,h}^{BSF}$\\\hline%
$\frac{\sqrt{2}}{4}$ & 33 & 2248.988  & 3401.184  & 20 & 8095.501  & 12267.987 \\
$\frac{\sqrt{2}}{8}$ & 161 & 4197.762  & 7026.946  & 132 & 7493.573  & 11830.283 \\
$\frac{\sqrt{2}}{16}$ & 705 & 5844.849  & 9933.029  & 644 & 7324.642  & 11784.194 \\
$\frac{\sqrt{2}}{32}$ & 2945 & 6680.697  & 11225.122  & 2820 & 7252.764  & 11777.697 \\
$\frac{\sqrt{2}}{64}$ & 12033 & 7000.825  & 11627.366  & 11780 & 7219.475  & 11776.312 \\
$\frac{\sqrt{2}}{128}$ & 48641 & 7116.044  & 11736.937  & 48132 & 7203.839  & 11775.948 \\
$\frac{\sqrt{2}}{256}$ & 195585 & 7159.231  & 11765.684  & ---- & ---- & ----\\
  \hline
\end{tabular}
\end{table}

\begin{table}
\small \centering {\bf Table 6: $\tau=100$, the 1st and the 2nd
eigenvalue on the L-shaped domain}\\
\begin{tabular}{ccccccc}\hline
h&Ndof&$\lambda_{1,h}^{M}$&$\lambda_{2,h}^{M}$&Ndof&$\lambda_{1,h}^{BSF}$&$\lambda_{2,h}^{BSF}$\\\hline%
$\frac{\sqrt{2}}{4}$ & 33 & 4138.743  & 6220.391  & 20 & 12775.204  & 19033.258 \\
$\frac{\sqrt{2}}{8}$ & 161 & 6676.425  & 10808.767  & 132 & 11854.535  & 18270.540 \\
$\frac{\sqrt{2}}{16}$ & 705 & 9245.708  & 15103.285  & 644 & 11635.005  & 18198.392 \\
$\frac{\sqrt{2}}{32}$ & 2945 & 10655.575  & 17239.122  & 2820 & 11548.155  & 18189.699 \\
$\frac{\sqrt{2}}{64}$ & 12033 & 11189.250  & 17931.933  & 11780 & 11508.527  & 18188.136 \\
$\frac{\sqrt{2}}{128}$ & 48641 & 11370.095  & 18121.348  & 48132 & 11489.938  & 18187.754 \\
$\frac{\sqrt{2}}{256}$ & 195585 & 11432.958  & 18170.592  & ---- & ---- & ----\\
  \hline
\end{tabular}
\end{table}

\begin{table}
\centering \small {\bf Table 7: $\tau=0, \theta=0.25$, the 1st and
the 2nd eigenvalue on the L-shaped domain}\\
\begin{tabular}{p{0.2cm}p{0.7cm}p{1.1cm}p{0.7cm}p{1.3cm}p{0.2cm}p{0.8cm}p{1.1cm}p{0.9cm}p{1cm}}
  \hline
 $l$ & Ndof & $\lambda_{1,h_{l}}^{M}$ & $Ndof$ & $\lambda_{2,h_{l}}^{M}$ &  $l$ & $Ndof$ & $\lambda_{1,h_{l}}^{M}$ & $Ndof$ & $\lambda_{2,h_{l}}^{M}$\\\hline
1 & 3201 & 6218.84  & 3201 & 10531.78  & 20 & 30033 & 6688.32  & 47187 & 11038.48 \\
2 & 3249 & 6345.81  & 3449 & 10661.95  & 21 & 35293 & 6690.89  & 54673 & 11041.15 \\
3 & 3403 & 6416.21  & 3801 & 10747.72  & 22 & 41310 & 6693.17  & 64375 & 11043.71 \\
4 & 3531 & 6456.17  & 4325 & 10828.75  & 23 & 48473 & 6695.17  & 75049 & 11045.12 \\
5 & 3829 & 6512.11  & 4950 & 10859.13  & 24 & 56871 & 6696.34  & 87881 & 11046.71 \\
6 & 4173 & 6536.95  & 5697 & 10847.20  & 25 & 66047 & 6697.07  & 103923 & 11047.31 \\
7 & 4643 & 6578.75  & 6543 & 10888.37  & 26 & 77216 & 6697.76  & 122683 & 11048.84 \\
8 & 5131 & 6584.41  & 7641 & 10923.11  & 27 & 90500 & 6698.80  & 142365 & 11049.90 \\
9 & 5793 & 6605.91  & 8905 & 10952.44  & 28 & 107015 & 6699.46  & 163283 & 11050.78 \\
10 & 6569 & 6623.88  & 10304 & 10966.49  & 29 & 126197 & 6700.25  & 187863 & 11051.26 \\
11 & 7623 & 6637.30  & 11917 & 10981.19  & 30 & 148263 & 6700.83  & 217569 & 11051.78 \\
12 & 8803 & 6649.11  & 13827 & 10987.69  & 31 & 174145 & 6701.28  & 255497 & 11052.19 \\
13 & 10201 & 6658.97  & 16175 & 10999.58  & 32 & 204237 & 6701.68  & 298477 & 11052.51 \\
14 & 11753 & 6661.76  & 18701 & 11008.80  & 33 & 240145 & 6701.95  & 348699 & 11052.76 \\
15 & 13621 & 6667.77  & 21889 & 11016.67  & 34 & 278865 & 6702.09  & 412401 & 11053.01 \\
16 & 15833 & 6674.39  & 25995 & 11022.76  & 35 & 326744 & 6702.30  & 487223 & 11053.31 \\
17 & 18437 & 6679.04  & 30840 & 11028.20  & 36 & 383793 & 6702.48  & 564663 & 11053.58 \\
18 & 21655 & 6683.12  & 35673 & 11031.47  & 37 & 454113 & 6702.68  & ----- & ----- \\
19 & 25453 & 6685.75  & 40981 & 11034.73  & 38 & 533383 & 6702.84  & ----- & ----- \\
  \hline
\end{tabular}
\end{table}

\begin{table}
\small \centering {\bf Table 8: $\tau=10,\theta=0.25$, the 1st and
the 2nd eigenvalue on the L-shaped domain}\\
\begin{tabular}{p{0.2cm}p{0.7cm}p{1.1cm}p{0.7cm}p{1.3cm}p{0.2cm}p{0.8cm}p{1.1cm}p{0.9cm}p{1cm}}
  \hline
 $l$ & Ndof & $\lambda_{1,h_{l}}^{M}$ & $Ndof$ & $\lambda_{2,h_{l}}^{M}$ &  $l$ & $Ndof$ & $\lambda_{1,h_{l}}^{M}$ & $Ndof$ & $\lambda_{2,h_{l}}^{M}$\\\hline
1 & 3201 & 6675.63  & 3201 & 11215.23  & 20 & 30733 & 7173.63  & 46779 & 11757.96 \\
2 & 3257 & 6809.63  & 3453 & 11356.72  & 21 & 36195 & 7176.68  & 54329 & 11761.11 \\
3 & 3423 & 6891.61  & 3817 & 11451.54  & 22 & 42307 & 7179.14  & 63757 & 11764.06 \\
4 & 3575 & 6934.33  & 4361 & 11532.55  & 23 & 49763 & 7181.08  & 74337 & 11765.37 \\
5 & 3897 & 6994.22  & 4985 & 11559.63  & 24 & 58131 & 7182.18  & 86857 & 11766.90 \\
6 & 4237 & 7013.90  & 5719 & 11551.52  & 25 & 67479 & 7182.84  & 102755 & 11767.89 \\
7 & 4735 & 7052.88  & 6569 & 11598.94  & 26 & 79069 & 7183.86  & 121155 & 11769.32 \\
8 & 5231 & 7067.93  & 7635 & 11630.88  & 27 & 92839 & 7184.81  & 140651 & 11770.75 \\
9 & 5919 & 7088.50  & 8845 & 11665.70  & 28 & 110135 & 7185.68  & 161305 & 11771.63 \\
10 & 6722 & 7104.83  & 10249 & 11678.69  & 29 & 129399 & 7186.42  & 185395 & 11772.28 \\
11 & 7755 & 7118.38  & 11829 & 11694.74  & 30 & 152321 & 7187.02  & 214803 & 11772.77 \\
12 & 9005 & 7132.88  & 13735 & 11701.14  & 31 & 178795 & 7187.51  & 251607 & 11773.13 \\
13 & 10417 & 7141.33  & 16063 & 11715.79  & 32 & 209879 & 7187.95  & 293607 & 11773.53 \\
14 & 11957 & 7145.09  & 18580 & 11725.40  & 33 & 245095 & 7188.13  & 342775 & 11773.81 \\
15 & 13906 & 7152.91  & 21749 & 11734.66  & 34 & 285181 & 7188.38  & 404877 & 11774.09 \\
16 & 16120 & 7158.40  & 25806 & 11740.72  & 35 & 334009 & 7188.54  & 478907 & 11774.47 \\
17 & 18869 & 7163.39  & 30511 & 11746.14  & 36 & 392889 & 7188.76  & 556009 & 11774.73 \\
18 & 22219 & 7168.25  & 35367 & 11749.89  & 37 & 464649 & 7188.96  & 638003 & 11774.89 \\
19 & 26145 & 7170.98  & 40623 & 11754.05  & 38 & 545321 & 7189.13  & ----- & ----- \\
\hline
\end{tabular}
\end{table}

\begin{table}
\small \centering {\bf Table 9: $\tau=100,\theta=0.25$,the 1st and
the 2nd eigenvalue on the L-shaped domain}\\
\begin{tabular}{p{0.2cm}p{0.7cm}p{1.3cm}p{0.7cm}p{1.3cm}p{0.2cm}p{0.8cm}p{1.3cm}p{0.9cm}p{1cm}}
  \hline
   $l$ & Ndof & $~~~\lambda_{1,h_{l}}^{M}$ & $Ndof$ & $~~~\lambda_{2,h_{l}}^{M}$ &  $l$ & $Ndof$ & $~~~\lambda_{1,h_{l}}^{M}$ & $Ndof$ & $~~~\lambda_{2,h_{l}}^{M}$\\\hline

1 & 3201 & 10647.46  & 3201 & 17225.37  & 20 & 35207 & 11445.82  & 46449 & 18153.22 \\
2 & 3297 & 10862.62  & 3461 & 17471.42  & 21 & 41155 & 11451.13  & 54005 & 18159.18 \\
3 & 3491 & 10982.42  & 3873 & 17627.06  & 22 & 47783 & 11454.61  & 62889 & 18162.96 \\
4 & 3739 & 11066.38  & 4452 & 17718.21  & 23 & 55501 & 11457.27  & 72625 & 18165.83 \\
5 & 4165 & 11139.87  & 5095 & 17736.72  & 24 & 64329 & 11458.63  & 84435 & 18168.31 \\
6 & 4613 & 11203.13  & 5807 & 17810.09  & 25 & 74704 & 11460.67  & 99655 & 18171.30 \\
7 & 5135 & 11230.75  & 6651 & 17873.26  & 26 & 87138 & 11462.50  & 116585 & 18174.04 \\
8 & 5782 & 11271.90  & 7711 & 17935.23  & 27 & 102486 & 11464.16  & 135147 & 18176.30 \\
9 & 6661 & 11305.70  & 8987 & 17976.08  & 28 & 120434 & 11465.59  & 155559 & 18178.18 \\
10 & 7673 & 11331.41  & 10500 & 18004.69  & 29 & 141483 & 11466.88  & 178829 & 18179.29 \\
11 & 8911 & 11355.96  & 12157 & 18024.75  & 30 & 165055 & 11467.85  & 206719 & 18180.98 \\
12 & 10391 & 11374.06  & 14036 & 18057.91  & 31 & 191875 & 11469.01  & 240853 & 18181.68 \\
13 & 12041 & 11386.79  & 16173 & 18079.99  & 32 & 222691 & 11469.55  & 278885 & 18182.69 \\
14 & 13923 & 11401.31  & 18775 & 18096.78  & 33 & 257621 & 11469.99  & 323925 & 18183.08 \\
15 & 16129 & 11412.36  & 21855 & 18107.55  & 34 & 298901 & 11470.38  & 380575 & 18183.78 \\
16 & 18701 & 11420.57  & 25857 & 18116.00  & 35 & 348221 & 11470.77  & 446081 & 18184.40 \\
17 & 21847 & 11428.60  & 30199 & 18126.82  & 36 & 408785 & 11471.15  & 517651 & 18184.90 \\
18 & 25721 & 11434.06  & 34921 & 18136.66  & 37 & 480267 & 11471.48  & 596345 & 18185.36 \\
19 & 30135 & 11438.79  & 40253 & 18146.53  & 38 & 563307 & 11471.86  & ----- & ----- \\
\hline
\end{tabular}
\end{table}

\begin{figure}
  \centering
  \includegraphics[width=6cm]{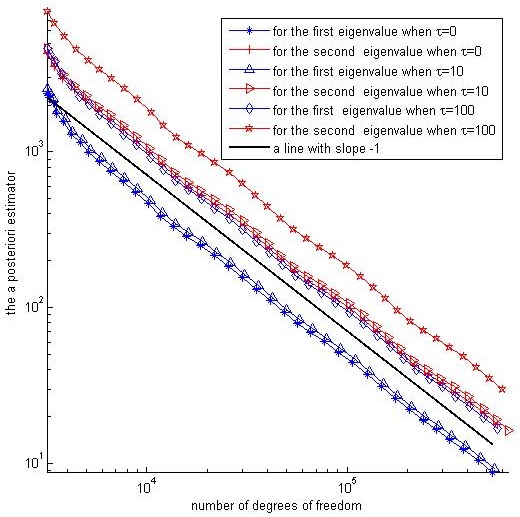}\\
  \caption{$\tau=0,10,100$, a posteriori error curve of the 1st and the 2nd
eigenvalue on the L-shaped domain}
\end{figure}

\section{Concluding remarks}

\indent In this paper, for fourth-order elliptic eigenvalue problems
with the clamped boundary condition in any dimension, including the
vibrations of a clamped plate under tension, we discuss the lower
bound property of the Morley element eigenvalues. We prove on
regular meshes, including adaptive local refined meshes, that in the
asymptotic sense the Morley element eigenvalues approximate the
exact ones from below, which is a development of the work in
\cite{hu1,lin2,yang2}. Our analysis is an application of the
identity (3.1), in which higher order
contributions are identified.\\
\indent Our analysis and results about the lower bound property of
the Morley element eigenvalues can be applied to general $2m$-th
order elliptic eigenvalue problems (see (6.5) in \cite{hu1}).

\indent{\bf Acknowledgements.}~~This work was supported by the National Natural Science Foundation of China (Grant Nos.11161012,11201093).\\


\footnotesize




\begin{thebibliography}{s37}
\bibitem{ainsworth}M. Ainsworth, J.T. Oden, A posteriori error
estimates in the finite element analysis. New York: Wiley-Inter
science, 2011.

\bibitem{ainsworth2}M. Ainsworth, Robust a posteriori error
estimation for nonconforming  finite element approximation, SIAM J.
Numer. Anal., 42 (2005): 2320-2341.

\bibitem{armentano}M.G. Armentano, R.G. Duran,  Asymptotic lower bounds
for eigenvalues by nonconforming finit element methods, Electron
Trans Numer. Anal., 17 (2004), pp. 92-101.


\bibitem{babuska}I. Babu$\check{s}$ka, J. Osborn, Eigenvalue Problems, Handbook of Numerical Analysis. vol.2, North-Holand:
Elsevier Science Publishers B.V., 1991.

\bibitem{babuska2}I. Babu$\check{s}$ka, R.B. Kellog, Pitkaranta J. Direct and inverse error estimates for finite elements
with mesh refinement, Numer. Math., 33 (1979), pp. 447-471.

\bibitem{babuska3}I. Babu$\check{s}$ka, A. Miller, A feedback finite element method with a posteriori error estimation:
Part I. The finite element method and some basic properties of the a
posteriori error estimator, Comput. Methods Appl. Mech. Engrg., 61
(1987), pp. 1-40.

\bibitem{bl}Veiga L. Beir$\tilde{a}$o da, J. Niiranen, R. Stenberg, A posteriori error estimates for the Morley plate
bending element, Numer. Math., 106 ( 2007), pp. 165-179.

\bibitem{becker2}R. Becker, S. Mao, Z. Shi, A Convergent nonconforming
adaptive finite element method with quasi-optimal complexity, SIAM
J. Numer. Anal., 47(6) (2010), pp. 4639-4659.

\bibitem{blum}H. Blum, R. Rannacher, On the boundary value problem of the biharmonic operator on domains with
angular corners, Math. Method. Appl. Sci., 2 (1980), pp. 556-581.

\bibitem{brenner}S.C. Brenner, L.R. Scott, The Mathematical Theory of Finite Element Methods. 2nd ed.,
Springer-Verlag, New york, 2002.

\bibitem{carstensen}C. Carstensen, D. Gallistl, Guaranteed lower eigenvalue bounds for the biharmonic
equation. Numer. Math., 1 (2014), pp. 33-51.

\bibitem{carstensen2}C. Carstensen, J. Hu, A. Orlando, Framework for
the a posteriori  error analysis of nonconforming finite elements.
SIAM J. Numer. Anal., 45(1) (2007), pp. 68-82.

\bibitem{chenl}L. Chen, iFEM: an innovative finite element methods package in MATLAB,
Technical Report, University of California at Irvine, 2009.

\bibitem{ciarlet}P.G. Ciarlet, Basic error estimates for elliptic proplems. Handbook of Numerical Analysis, vol.2, North-Holand:
Elsevier Science Publishers B.V., 1991.

\bibitem{dai1}X. Dai, J. Xu, A. Zhou, Convergence and optimal complexity of
adaptive finite element eigenvalue computations. Numer. Math., 110
(2008), pp. 313-355.

\bibitem{duran}R.G. Dur$\acute{a}$n, C. Padra, R. Rodr$\acute{i}$guez, A posteriori error estimates
for the finite element approximation of eigenvalue problems, Math.
Mod. Meth. Appl. Sci., 13 (2003), pp.1219-1229.

\bibitem{gallistl}D. Gallistl, Morley Finite Element Method for the Eigenvalues of the
Biharmonic Operator, arXiv:1406.2876 [math.NA], 2014

\bibitem{gedicke} J. Gedicke and C. Carstensen, A posteriori error
estimators for convection-diffusion eigenvalue problems, Comput.
Methods Appl. Mech. Engrg., 268(2014), pp.160-177.

\bibitem{hu1}J. Hu, Y. Huang, Q. Lin, Lower bounds for eigenvalues of elliptic
operators: By Nonconforming finite element methods, J. Sci. Comput.,
61(2014), pp.196-221

\bibitem{hu2}J. Hu, Z. Shi, A new a posteriori error estimate for the Morley element, Numer. Math.,
112(1)(2009), pp. 25-40.

\bibitem{lascaux}P. Lascaux, P. Lesaint, Some nonconforming finite elements for the plate bending problem,
Rev. Francaise Automat. Informat. Recherche Operationnelle Ser Rouge
Anal Numer, R-1(1975), pp. 9-53.

\bibitem{lin2}Q. Lin, H. Xie, Recent results on lower bounds of
eigenvalue problems by nonconforming finite element methods. Inverse
Problems and Imaging, 7(3) (2013), pp. 795-811.

\bibitem{lin3}Q. Lin, H. Xie, J. Xu, Lower bounds of the discretization for piecewise polynomials,
Math. Comp., 83(285) (2014), pp. 1-13.

\bibitem{luo}F. Luo, Q. Lin, H. Xie, Computing the lower and upper bounds of
Laplace eigenvalue problem: by combining conforming and
nonconforming finite element methods. Sci. China Math., 55 (2012),
pp. 1069-1082.

\bibitem{mao}S. Mao, S. Nicaise, Z. Shi, Error estimates of Morley triangular element
satisfying the maximal angle condition. Int. J. Numer. Anal. Model.,
7 (2010), pp. 639-655.

\bibitem{morin}P. Morin, R.H. Nochetto, K. Siebert, Convergence of
adaptive finite element methods. SIAM Rev, 44 (2002), pp. 631-658.

\bibitem{morley} L.S.D. Morley, The triangular equilibrium element in the solution of plate bending problems, Aero Quart
19 (1968), pp. 149-169.

\bibitem{oden}J. T. Oden, J. N. Reddy, An Introduction to the Mathematical Theory of Finite Elements,
Courier Dover Publications, New York, 2012.

\bibitem{rr2}R. Rannacher, Nonconforming finite element methods for eigenvalue
problems in linear plate theory, Numer. Math., 33 (1979), pp. 23-42.

\bibitem{shen}Q. Shen A posteriori error estimates of the Morley element for the
fourth order elliptic eigenvalue problem, Numer. Algor., 68 (2015),
pp. 455-466

\bibitem{shi1}Z. Shi, On the error estimates of Morley element. Chinese J. Numer. Math. \& Appl., 12 (1990), pp.
102-108.

\bibitem{shi2}Z. Shi, M. Wang, Finite Element Methods. Beijing,
Scientific Publishers, 2013.

\bibitem{strang}G. Strang, G.J. Fix, An alalysis of the finite element method. Prentice-Hall,
1973.

\bibitem{verfurth}R. Verf$\ddot{u}$rth, A review of a posteriori
error estimates and adaptive mesh-refinement techniques. New York:
Wiley-Teubner, 1996.

\bibitem{wm}M. Wang, J. Xu, The Morley element for fourth order elliptic equations in any
dimensions. Numer. Math., 103(2006), pp. 155-169.

\bibitem{wm2}M. Wang, J. Xu, Minimal finite-element spaces for 2m-th order partial differential equations in $R^{n}$.
 Math Camp, 82(2013), pp. 25-43.

\bibitem{weinstein}A. Weinstein, W.Z. Chien, On the vibrations of a clamped plate under
tension, Quarterly of Applied Mathematics, 1(1) (1943), pp. 61-68.

\bibitem{widlund}O. Widlund, On best error bounds for approximation by
piecewise polynomial functions. Numer. Math., 27 (1977), pp.
327-338.

\bibitem{yang1}Y. Yang, Z. Zhang, F. Lin, Eigenvalue approximation from below using
nonforming finite elements. Sci. China Math., 53(1) (2010), pp.
137-150.

\bibitem{yang2}Y. Yang, Q. Lin, H. Bi, Q. Li, Lower eigenvalues approximation by Morley
elements. Adv. Comput. Math., 36(3) (2012), pp. 443-450.

\bibitem{yang4}Y. Yang, J. Han, H. Bi, Y. Yu, The lower/upper bound
property of the nonoconforming Crouzeix-Raviart element eigenvalues
on adaptive meshes. J. Sci. Comput., 62(1) (2015), pp. 284-299.

\bibitem{yang5}Y. Yang, A posteriori error analysis of
conforming/nonconforming finite elements (in Chinese), Sci. Sin.
Math., 40(9) (2010), pp. 843-862.

\bibitem{zhang}Z. Zhang, Y. Yang, Z. Chen, Eigenvalue approximation from below
by Wilson's element. Chinese J. Numer. Math. \& Appl., 29(4) (2007),
pp. 81-84.



\end{thebibliography}
\end{document}